\documentclass{article}
\usepackage[left=3.6cm,right=3.6cm]{geometry}
\usepackage{amssymb,amsthm,amsmath,amsxtra,dsfont,appendix}
\usepackage[nottoc]{tocbibind}   
\usepackage{mathrsfs}

\numberwithin{equation}{section} 

\renewcommand{\limsup}[1]{\underset{{#1}}{\overline{\operatorname{lim}}}}
\renewcommand{\liminf}[1]{\underset{{#1}}{\underline{\operatorname{lim}}}}

\newcommand{\sgn}{\operatorname{sgn}}
\newcommand{\supp}{\operatorname{supp}}
\newcommand{\tr}{\operatorname{Tr}}
\newcommand{\Var}{\operatorname{Var}}

\newcommand{\E}[1]{\mathbb{E}\left[#1\right]}

\newcommand{\1}{\mathds{1}}
\newcommand{\C}{\mathbb{C}}
\newcommand{\Cu}{\operatorname{C}}

\newcommand{\G}{\operatorname{GUE}}
\newcommand{\N}{\mathbb{N}}
\newcommand{\No}{\mathcal{N}}
\newcommand{\R}{\mathbb{R}}
\newcommand{\Sy}{\mathbb{S}}
\newcommand{\s}{\operatorname{sc}} 
\newcommand{\Z}{\mathrm{Z}}

\newtheorem{theorem}{Theorem}[section]
\newtheorem{definition}[theorem]{Definition}
\newtheorem{proposition}[theorem]{Proposition}

\newtheorem{lemma}[theorem]{Lemma}

\theoremstyle{definition} \newtheorem{remark}{Remark}[section]

\title{ \vspace{2cm} 
Mesoscopic fluctuations for unitary invariant
ensembles}
\date{}
\author{Gaultier Lambert\thanks{KTH Royal Institute of Technology, Department of Mathematics, glambert@kth.se. Supported by the grant KAW 2010.0063 from the Knut and Alice Wallenberg Foundation.}}

\begin{document}

\maketitle

\vspace{2cm}

\begin{abstract} \normalsize
\noindent
Considering a determinantal point process on the real line, we establish a connection between the sine-kernel asymptotics for the correlation kernel and the CLT for mesoscopic linear statistics. This implies universality of mesoscopic fluctuations for a large class of unitary invariant Hermitian ensembles.
In particular, this shows that the support of the equilibrium measure need not be connected in order to see Gaussian fluctuations  at mesoscopic scales.    
Our proof is based on the cumulants computations introduced in~\cite{Soshnikov_00a} for the CUE and the sine process and the asymptotics formulae derived by Deift et al.~\cite{Deift_al_99_b}. 
For varying weights $e^{-N \tr  V (H)}$, in the one-cut regime, we also provide estimates for the variance of linear statistics $\tr f(H)$ which are valid for a rather general function~$f$. In particular, this implies that the characteristic polynomials of such Hermitian random matrices converge in a suitable regime to a regularized fractional Brownian motion with logarithmic correlations introduced in~\cite{FKS_13}.  
For the GUE kernel, we also discuss how to obtain the necessary sine-kernel asymptotics at mesoscopic scale by elementary means. \\
\end{abstract}

\vspace{1.5cm}



\clearpage

\
\vspace{1cm}

\tableofcontents

\clearpage

\section{Introduction and results}


A point process on $\R$ is called determinantal if its correlation functions (with respect to the Lebesgue measure) exist and are given by  
\begin{equation} \label{det} \rho_k(x_1,\dots,x_k) =  \det_{k\times k} \big[K(x_i,x_j)\big]  \ ,\hspace{.8cm} \forall x_1,\dots,x_k\in\R \ ,  \hspace{.8cm}  \forall  k\in\N \ , 
\end{equation}
where $K:\R\times\R \to \R$ is called the correlation kernel.
These processes arise in random matrix theory  to describe eigenvalues of the so-called unitary (invariant) ensembles; see theorem~\ref{thm:UIE} below and section~\ref{sect:RMT} for more details. There are many other interesting examples such as random tilings or the positions of non-colliding stochastic diffusions that we will not discuss here. We refer to ~\cite{Soshnikov_00b,Johansson_05,HKPV_06,Borodin_11,Lyons_14} for various introductions to the general theory and further examples. A fundamental feature of determinantal processes is that all the information about the process is encoded in the correlation kernel. For instance, for unitary invariant ensembles, universality of the local correlations in the bulk of the spectrum follows from the convergence of the rescaled correlation kernel to the sine-kernel, \cite{Deift_99, Pastur_Shcherbina, SV_15}. In this work, we show that at mesoscopic scales, the sine-kernel  asymptotics still holds  and this leads to  the following Central Limit Theorem.


\begin{theorem}\label{thm:UIE}   Let $V:\R\to\R$ be real-analytic such that 
\begin{equation} \label{potential} \liminf{|x|\to\infty} \frac{V(x)}{\log(x^2+1)} = +\infty \ ,
\end{equation}
and consider the probability measure  $d\mathbb{P}^V_N =  \Z_{V,N}^{-1}e^{-N \tr  V (H)} dH  $ on the space of $N\times N$ Hermitian  matrices equipped with the Lebesgue measure $dH$. If $(\lambda_1, \dots, \lambda_N)$ denote  the eigenvalues of a random matrix $H$ distributed according to  $\mathbb{P}^V_N$, then for any $x_0\in J_V$, any  $0<\alpha<1$, and for any $f\in C^1(\R)$ with compact support, we have as $N\to\infty$, 
\begin{equation} \label{clt_1}
 \sum_{k=1}^N f\big( N^{\alpha}(\lambda_k-x_0)\big) - \mathbb{E}^V_N\left[  \sum_{k=1}^N f\big( N^{\alpha}(\lambda_k-x_0)\big) \right] \ \Rightarrow \ \No\big( 0, \|f\|^2_{H^{1/2}} \big) \ .
\end{equation}
\end{theorem}

\proof Section~\ref{sect:UIE}.\qed\\

The condition (\ref{potential}) guarantees that  $\Z_{V,N} <\infty$, so that the measure $\mathbb{P}^V_N$ is well-defined. This also implies that, for large $N$, the eigenvalue process  is supported in a fixed compact set $\overline{J_V}$ with high probability; see formula (\ref{equilibrium}) below. Hence, the potential $V$ is confining and the condition $x_0\in J_V$  means that we zoom in around a point $x_0$ which lies in the bulk of the spectrum. \\


In theorem \ref{thm:UIE}, the parameter $\alpha\in[0,1]$ is called the scale. Since the eigenvalues density is of order $N$ in the bulk, when $\alpha=0$,
  the l.h.s$.$ of (\ref{clt_1}) depends on the whole spectrum of $H$ and this regime is called  global or macroscopic.
On the other hand, when $\alpha=1$, the distance between consecutive eigenvalues remains of order 1 as  the size $N$ of the matrix tends to infinity and this regime is called local or microscopic.  Any intermediate scale, $0<\alpha<1$, is called {\bf mesoscopic}. In this regime, the limit (\ref{clt_1}) is independent of the potential $V$, the scale $\alpha$ and $x_0$. Hence, this establishes universality of fluctuations for a large class of Hermitian random matrix ensembles.


The variance in formula (\ref{clt_1}) is given by
\begin{equation}\label{norm_1}
\|f\|^2_{H^{1/2}} = \int_{\R} \big|\hat f (u) \big|^2 |u| du =  \frac{1}{4\pi^2}  \iint_{\R^2} \left| \frac{f(x)- f(y)}{x-y} \right|^2 dxdy \ , 
\end{equation}
where $\displaystyle \hat f(u)= \int f(x) e^{-i2\pi x u} dx $ denotes the Fourier transform of $f$.
Modulo constant functions, the norm (\ref{norm_1}) defines a complete subspace of $L^2(\R)$ denoted by $H^{1/2}(\R)$.  
Most of the work on unitary invariant ensembles has focused on the asymptotics of local or global statistics and we briefly review the main results, further references can be found in the textbooks~\cite{Deift_99, Pastur_Shcherbina}.
Under the assumptions of theorem~\ref{thm:UIE}, there exists a probability density $\varrho_V$ with compact support $\overline{J_V}$ on $\R$ such that for any $f\in C\cap L^\infty(\R)$,
\begin{equation} \label{equilibrium}
 \frac{1}{N} \sum_{k=1}^N f(\lambda_k) \underset{N\to\infty}{\longrightarrow}  \int f(x) \varrho_V(x) dx    \ , \hspace{1cm}  \mathbb{P}^V_N \text{ - almost surely.}
\end{equation}

It means that, for large $N$, the eigenvalues of a random matrix sampled according to $\mathbb{P}^V_N$ are distributed according to the {\bf equilibrium density} $\varrho_V$. 
Moreover, it is known that the fluctuations around this equilibrium configuration remain bounded as $N\to\infty$. 
The precise behavior of linear statistics depends on the support of $\varrho_V$. In the simplest case, there exists $x_0\in\R$ and  $ \ell>0$ so that
\begin{equation}\label{1cut}
J_V = \big( x_0 - \ell, x_0 + \ell) \ , 
\end{equation}
the potential $V$ is said to satisfy the {\bf one-cut} condition and we have a CLT: for  any $f\in C^2\cap L^\infty(\R)$,
\begin{equation} \label{clt_2}
\sum_{k=1}^N f\big(\frac{\lambda_k-x_0}{\ell}\big)  - \ell \int f(x) \varrho_V\big(x_0+ \ell x\big) dx  \underset{N\to\infty}{\Longrightarrow}   \No\big( 0, \Sigma(f)^2 \big)    \ , 
\end{equation}
where 
\begin{equation}\label{norm_2} 
\Sigma(f)^2:= \frac{1}{4\pi^2}  \iint\limits_{[-1,1]^2} \left| \frac{f(x)- f(y)}{x-y} \right|^2 \frac{1- xy}{\sqrt{1-x^2}\sqrt{1-y^2}} dxdy  \ .
\end{equation} 

Theorem (\ref{clt_2}) was first proved in~\cite{Johansson_98} when $V$ is a polynomial of even degree using a variational method. We refer to~\cite{BG_13a} for further developments and to \cite{BD_13, L_15b} for alternative proofs which are valid for more general determinantal processes. 
It is known that (\ref{1cut}) holds when the potential $V$ is strictly convex on $\R$ and, if $\tilde{V}(x)=V(x_0+\ell x)$, by considering the ensemble $\mathbb{P}^{\tilde V}_N$ instead of $\mathbb{P}^V_N$, we can always assume that $x_0=0$ and $\ell=1$.
The one-cut condition is crucial to observe a Gaussian process in the limit. If $\overline{J_V}=\supp(\varrho_V)$ is not connected, then  for a generic test function~$f$,  the behavior of the linear statistic $\sum f(\lambda_k)$ is quasi-periodic in $N$ and, even though this sequence of random variables is tight, it has no limit has $N\to\infty$, see~\cite{Pastur_06}. This complicated behavior is explained by the fact that the number of eigenvalues in the different components of $J_V$ fluctuates. Nevertheless, along the subsequence where it converges, the  asymptotic distribution of $\sum f(\lambda_k)$  can still be described and it is not Gaussian in general, \cite{BG_13b}. 
On the other hand, at the local scale, the behavior of the eigenvalue process  is independent of the equilibrium density and it is described by the sine process in the bulk. Theorem~\ref{thm:UIE} shows that mesoscopic fluctuations are universal as well. Actually, this results was first derived heuristically by Pastur in~\cite{Pastur_06} also based on the {\it semiclassical} asymptotic formulae derided in~\cite{Deift_al_99_b} for the orthogonal polynomials with respect to the measure $e^{-N V(x)}dx$ on~$\R$.     \\

 Mesoscopic spectral statistics  were first considered in \cite{BK_99a,BK_99b}  for Hermitian and symmetric Wigner matrices.  In particular, the authors proved a result analogous to (\ref{clt_1}) for the GUE and GOE using the resolvent as a test function.  
One of the pioneering works on the subject which has been of inspiration for this article is  Soshnikov's CLT for eigenvalues statistics of Haar distributed random matrices from the compact groups, \cite{Soshnikov_00a}. In the case of the unitary group, this probability measure is known as the Circular Unitary Ensemble (CUE) and it is determinantal with the correlation kernel
\begin{equation}\label{CUE}
 K_{U_N}(x,y)=\frac{\sin \big( N (x-y)/2 \big)}{2 \pi \sin\big((x-y)/2\big)}  \hspace{1cm} \forall x, y \in \R/ 2\pi\mathbb{Z} \ .  
\end{equation}   
For this point process, Soshnikov obtained the counterpart of (\ref{clt_1}) which can been seen as a continuous analogue of the Strong Szeg\H{o} theorem.
 A special case of  theorem~\ref{thm:UIE} was recently established in \cite{FKS_13, LS_15} for the Gaussian Unitary Ensemble (GUE). The authors of~\cite{FKS_13} proved that a suitable regularization of the  characteristic polynomial of a GUE matrix  converges weakly to a certain {\it fractional Brownian motion} which is logarithmically correlated, see section~\ref{sect:fbm}. From Theorem~2.4 therein, one  can infer the  CLT for mesoscopic linear statistics of any Schwartz-class  test function. In~\cite{LS_15}, analogous results are proved for  Hermitian Wigner matrices with sub-Gaussian entries, extending the results of \cite{BK_99b}. 
For a class of determinantal processes known as orthogonal polynomial ensembles,
 an alternative approach to universality  which is discussed in section~\ref{sect:intro} appeared in \cite{BD_13, BD_15}. In particular, the authors obtained the counterpart of  (\ref{clt_1}) for another family of unitary invariant ensembles, see theorem~\ref{thm:JE} and remark~\ref{rk:BD}. All these results have the following interpretation. If we view the eigenvalue process  as a random measure
\begin{equation} \label{pp}
 \Xi_N^{x_0,\alpha} := \sum_{k=1}^N \delta_{ N^{\alpha}(\lambda_k-x_0)} \ ,
\end{equation}
if centered, $ \Xi_N^{x_0,\alpha}$ converges in distribution to a random Gaussian process $\mathfrak{G}$  with covariance structure
\begin{equation}\label{noise}
\E{ \mathfrak{G}(f) \mathfrak{G}(g)} =  \int_{\R} \hat f (u) \hat g (-u)   |u| du \ .
\end{equation}

The random process $\mathfrak{G}$ is called the $H^{1/2}$-Gaussian field, see \cite{Janson_97}  chap.~1. Its special feature is that it is scale invariant. If $f_\eta(x)= f(\eta x)$, then $\mathfrak{G}( f_\eta) \sim \mathfrak{G}( f) $   for any $\eta>0$,   as can be seen from formula (\ref{noise}). 
Heuristically, this motivates why it is expected to describe mesoscopic fluctuations of point processes with strong repulsion such as  eigenvalues of random matrices, see the discussion in \cite{Soshnikov_01}. In some respect, these ensembles behave like the sine process and this is the idea behind the proof of  theorem \ref{thm:universality} below. The mesoscopic correlations can also be guessed from formulae (\ref{clt_2} - \ref{norm_2}). Namely, if $x_0=0$ and $\ell=1$, by a change of variables
\begin{equation*}
\Sigma(f_\eta)^2= \frac{1}{4\pi^2}  \iint\limits_{[-\eta,\eta]^2} \left| \frac{f(x)- f(y)}{x-y} \right|^2 \frac{1- xy/\eta^2}{\sqrt{1-(x/\eta)^2}\sqrt{1-(y/\eta)^2}} dxdy  \ ,
\end{equation*} 
and, if $f$ decays sufficiently fast,  we obtain
 \begin{equation*}
\lim_{\eta\to\infty}\Sigma(f_\eta)^2= \frac{1}{4\pi^2}  \iint\limits_{\R^2} \left| \frac{f(x)- f(y)}{x-y} \right|^2 dxdy
= \|f\|_{H^{1/2}}^2  \ .
\end{equation*} 


It is natural to investigate whether (\ref{clt_1}) holds under the optimal condition $f\in H^{1/2}\cap L^1(\R)$. To the author' knowledge, the question    remains open even for the Gaussian Unitary Ensemble (GUE).
To some extent, this issue is addressed in section~\ref{sect:variance}.
 In particular, if the potential $V$ satisfies the one-cut condition (\ref{1cut}), we derive an upper-bound for the variance of the random variable $\Xi_N^{x_0,\alpha}f$ which is valid  e.g.~for any function $f\in H^{1/2}(\R)$ with compact support, cf.~proposition~\ref{thm:variance_meso}.
  This easily allows us to generalize  theorem~\ref{thm:UIE}, cf.~theorem~\ref{thm:clt_1}.
 As a corollary,  we establish in section~\ref{sect:fbm}  that  the result of~\cite{FKS_13} is also valid for the characteristic polynomial of a random matrix from an arbitrary  unitary invariant  $\mathbb{P}^V_N$ in the one-cut regime.\\

The proof of theorem~\ref{thm:UIE} is based on the  so-called Plancherel-Rotach asymptotics for the orthogonal polynomials (OP) with respect to the weight $e^{-N V (x)}$ on $\R$  derived in~\cite{Deift_al_99_b} and the following general result.
%
%
For any function $\rho:\R\to\R_+$ locally integrable, we let 
\begin{equation} \label{bulk}
J_\rho := \left\{ t\in\R : \rho(t) >0 \ \text{and}\ \rho(t)\ \text{is continuous}\right\}  \,
\end{equation}
and, for all $x\in\R$, we define
\begin{equation} \label{F}
F_\rho(x): = \int_0^x \rho(t) dt \ .
\end{equation}
We also denote by $C^k_0(\R)$ the space of compactly supported real-valued functions with $k$ continuous derivatives on  $\R$.

\begin{theorem}\label{thm:universality}Consider a determinantal process on $\R$ with a correlation kernel $K_N$ which is locally trace-class. 
For any $x_0 \in \R$, $\alpha\in [0,1]$ and $f\in C_0(\R)$, we consider the linear statistic
\begin{equation*}
 \Xi_N^{x_0,\alpha} f  = \sum f\big( N^{\alpha}(\lambda_k-x_0)\big) \ ,
\end{equation*}
where the sum is over the point configuration $\{\lambda_k\}$ of the process.  
Assume that there exists a  function $\rho:\R\to\R_+$,  $x_0 \in J_\rho$,  $\alpha\in (0,1]$ and $\beta>0$ such that for any $L>0$,
\begin{equation} \label{sine}
\frac{1}{N^{\alpha}} K_N\big(x_0 + \frac{\xi}{N^{\alpha}} , x_0 + \frac{\zeta}{N^{\alpha}} \big)= \frac{\sin \pi N \big( ( F_\rho (x_0 + \xi N^{-\alpha} ) - F_\rho(x_0 + \zeta N^{-\alpha} )\big)}{\pi(\xi-\zeta)} + \underset{N\to\infty}{O}(N^{-\beta}) \ . 
\end{equation}
uniformly for all $\xi,\zeta \in [-L,L]$. Then, if $\alpha<1$, for any $f\in C^1_0(\R)$, we have as $N\to\infty$,
\begin{equation} \label{clt} \Xi_N^{x_0,\alpha} f - \E{ \Xi_N^{x_0,\alpha} f } \Rightarrow \No\left(0, \|f\|_{H^{1/2}}^2 \right) \ . 
\end{equation}
On the other hand, if $\alpha=1$,  for any $f\in C_0(\R)$, we have as $N\to\infty$,
\begin{equation} \label{slt}  \Xi_N^{x_0,1} f  \Rightarrow  \Xi^{\sin}_{\rho(x_0)} f \  .  
\end{equation}
\end{theorem}

\proof Section~\ref{sect:universality}. \qed\\

For any $\nu>0$, $  \Xi^{\sin}_\nu$ denotes the sine process with density $\nu>0$, i.e.~the determinantal process on $\R$ with the correlation kernel 
  \begin{equation} \label{sine_K}
 K^{\sin}_\nu(\xi,\zeta)=  \frac{\sin[\pi \nu(\xi-\zeta)]}{\pi(\xi-\zeta)} \ .
 \end{equation}
 
At the local scale, the limit (\ref{slt}) implies the convergence  of the process $ \Xi_N^{x_0,1}$ to the sine process. This behavior is known to be universal for Hermitian ensembles. In the context of theorem~\ref{thm:UIE}, it was proved in \cite{PS_97,Deift_al_99_b,LL_08,LL_09}. 
%
%
Assuming that the kernel $K_N$ is  locally trace-class is standard, it means that for any function $f\in L^\infty(\R)$ with compact support, the integral operator    
 $$ h \mapsto \int h(x) f(x) K_N(\cdot,x) dx  $$
is trace-class on $L^2(\R)$. For instance, this implies that the linear statistic $\sum f(\lambda_k)$ has a finite Laplace transform and its cumulants are well-defined, see formula (\ref{cumulant_0}) below.
Note that, we do not assume that the  kernel $K_N$ is  reproducing. In particular, the configuration $\{\lambda_k\}$ may have a random number of points or infinitely many. Hence, theorem~\ref{thm:universality}  can be applied beyond the context of unitary ensembles.\\     

It is obvious that the CUE kernel (\ref{CUE})  has an asymptotic expansion of the form (\ref{sine}) with $\rho=1/2\pi$ and Soshnikov's CLT is a special case of theorem~\ref{thm:universality}. 
Our main observation is that, if the correlation kernel satisfies (\ref{sine}), then we can still apply Soshnikov's method to prove a Central limit theorem, see lemma \ref{thm:cumulant}. In particular, the fact that the limiting process is Gaussian follows from the Main combinatorial Lemma of  \cite{Soshnikov_00a}, theorem~\ref{thm:MCL}.  
For determinantal processes within the sine process universality class, it is plausible that the asymptotics (\ref{sine}) holds at scales $\alpha$ which are sufficiently close to 1, so that theorem \ref{thm:universality} explains the appearance of the $H^{1/2}$-Gaussian field $\mathfrak{G}$ in this context. This also emphasizes on the connection with the Main combinatorial Lemma.  
However, the general mechanism behind universality of mesoscopic fluctuations is still far from being understood. In particular, it would be interesting to understand further the connection between random matrix theory and logarithmically correlated Gaussian fields as discussed in \cite{FKS_13, LS_15}. 
Within other symmetry classes and for the Dyson's $\beta$-ensembles, mesoscopic correlations are also conjectured to be described by the $H^{1/2}$-Gaussian field. For instance, this has been rigorously established for the Gaussian $\beta$-Ensembles in \cite{BEYY_14}, for random matrices from the special orthogonal and symplectic groups in \cite{Soshnikov_00a}  and in number theory, when considering mesoscopic linear statistics of the zeros of the Riemann-Zeta function \cite{BK_14,Rodgers_14}. 
There are also examples of determinantal processes where the precise asymptotics of the correlation kernel is  not known, but the CLT (\ref{clt}) has been proved by other means, e.g$.$ for non-colliding Brownian motions in \cite{DJ_14}.  \\

In this article, we focus on applications to random matrices, but  theorem~\ref{thm:universality} should be useful to investigate mesoscopic fluctuations for other instances of determinantal processes. Based on the 
Riemann-Hilbert formulation of \cite{FIK_92,  DZ_93}, it is possible to derive very precise asymptotics for  the orthogonal polynomials and the Christoffel-Darboux kernels for a large class of measures on $\R$. These results combined with theorem~\ref{thm:universality} allows to prove universality of the mesoscopic correlations for the classical random matrix ensembles.  For the GUE, it is possible to derive the asymptotics (\ref{sine}) using  only the Plancherel-Rotach asymptotics for the Hermite polynomials, \cite{PR_29}, and this leads to a rather elementary proof of theorem~\ref{thm:UIE}.

\clearpage 

The rest of the paper is organized as follows. In section~\ref{sect:universality}, we review the cumulants method introduced in \cite{Soshnikov_00a} to study linear statistics of determinantal processes and we prove theorem~\ref{thm:universality}. The proof relies on ideas developed in~\cite{JL_15}.
In section~\ref{sect:intro}, we give a brief introduction to the theory of unitary ensembles focusing on the orthogonal polynomials method. In section~\ref{sect:variance}, we provide some estimates for the variance of linear statistics of orthogonal polynomial ensembles both in the global and mesoscopic regime. 
 In section~\ref{sect:UIE}, we review  the results of \cite{Deift_al_99_b}  on the asymptotics of the Christoffel-Darboux kernels for varying exponential weights. This provides the necessary  asymptotics to prove theorem~\ref{thm:UIE}.
 In section~\ref{sect:JE}, we discuss another family of unitary ensembles known as the Moified Jacobi ensembles. The asymptotics of their correlation kernels is derived using the results of \cite{Kuijlaars_al_04, KV_02} and we deduce a CLT in this case as well, theorem~\ref{thm:JE}.
 In section~\ref{sect:fbm}, we apply the mesoscopic CLT to generalize the result of~\cite{FKS_13} to a large class of unitary invariant ensembles in the one-cut regime.   In section~\ref{sect:GUE}, we present an elementary approach to obtain the sine-kernel asymptotics at mesoscopic scales for the GUE kernel; see theorem~\ref{thm:sine_GUE}. In the appendix~\ref{A:variance}, we generalize an estimate obtained in section~\ref{sect:variance} for the variance of global linear statistics, further motivations and applications are given in~\cite{L_15b}.

\section{Proof of theorem~\ref{thm:universality}} \label{sect:universality}


We consider a family of determinantal processes on $\R$ with correlation kernels $K_N$ which depend on a parameter $N>0$.
 We want to study the law of the random variable
\begin{equation} \label{ls}
 \Xi_N^{x_0,\alpha} f  = \sum f\big( N^{\alpha}(\lambda_k-x_0)\big) \ ,
\end{equation}
 as $N\to\infty$, where $\{\lambda_k\}$ is a  configuration of the determinantal process with kernel $K_N$ and $f \in C^1_0(\R)$ is a test function. We will assume that $\supp(f) \subset [-L,L]$.\\

 For any real-valued random variable $Z$ with a well-defined Laplace transform, its cumulants $\Cu^n[Z]$ are defined by the generating function:  
\begin{equation} \label{cumulant}
\log \E{e^{tZ}} = \sum_{n=1}^\infty \Cu^n[Z] \frac{t^n}{n!}  \ .
\end{equation}     
Observe that $\Cu^1[Z]=\E{Z}$ and the higher-order cumulants do not depend on $\E{Z}$. In particular, we have $\Cu^2[Z]=\Var[Z]$ and, if $Z$ is Gaussian, $\Cu^n[Z]=0$ for all $n\ge3$. Hence, to prove the CLT (\ref{clt}), it is enough to show that
\begin{equation} \label{cumulant_limit}
 \lim_{N\to\infty} \Cu^2\big[ \Xi_N^{x_0,\alpha} f  \big] =   \|f\|_{H^{1/2}}^2  \ , \hspace{1cm}
 \lim_{N\to\infty} \Cu^n\big[ \Xi_N^{x_0,\alpha} f  \big] =  0  \hspace{.5cm} \forall n\ge3 \ .
\end{equation} 
 Using formula (\ref{det}), one can compute moments and cumulants of linear statistics of determinantal processes. In particular, it was proved in \cite{Soshnikov_00a} that, if the correlation kernel is locally trace-class and $f\in C_0(\R)$, then for any $n\in\N$,  
\begin{equation} \label{cumulant_0}
\Cu^n\left[ \sum f(\lambda_k)\right]
= \sum_{\ell=1}^n \frac{(-1)^{\ell+1}}{\ell} \sum_{\begin{subarray}{c} m_1, \dots, m_\ell \ge 1 \\ m_1+ \cdots + m_\ell = n  \end{subarray}  } \frac{n!}{m_1!\cdots m_\ell!}
  \tr\left[f^{m_1}K_N \cdots f^{m_\ell} K_N\right]  \ ,
  \end{equation}
where
\begin{equation}\label{trace}
  \tr\left[f^{m_1}K \cdots f^{m_\ell} K\right] 
  = \int_{\R^\ell} f(x_1)^{m_1}K(x_1,x_2)\cdots f(x_\ell)^{m_\ell} K(x_\ell,x_1) d^\ell x \ .
\end{equation}

\vspace{.3cm}

In the following, we suppose that $K_N$ satisfies (\ref{sine}) for a given function $\rho:\R\to\R^+$ and we let $J=J_\rho$ and $F=F_\rho$ according to (\ref{bulk}), respectively (\ref{F}).
We also assume that $J$ is non-empty, fix a point $x_0\in J$ and, for any $\xi,\zeta\in \R$, we let
 \begin{equation*} 
 \tilde K_N(\xi,\zeta) =  N^{-\alpha} K_N\left(x_0 +N^{-\alpha} \xi, x_0+N^{-\alpha} \zeta\right)  \ .
 \end{equation*} 
Then, by (\ref{ls}), (\ref{cumulant_0}) and a change of variables, we get 
\begin{equation} \label{cumulant_1}
\Cu^n\big[  \Xi_N^{x_0,\alpha} f  \big]
= \sum_{\ell=1}^n \frac{(-1)^{\ell+1}}{\ell} \sum_{\begin{subarray}{c} m_1, \dots, m_\ell \ge 1 \\ m_1+ \cdots + m_\ell = n  \end{subarray}  } \frac{n!}{m_1!\cdots m_\ell!}
  \tr\left[f^{m_1} \tilde K_N \cdots f^{m_\ell} \tilde K_N \right]  \ .
  \end{equation}

  
It was observed in  \cite{JL_15} that, if the correlation kernel $K_N$ satisfy the uniform asymptotics (\ref{sine}), then we can relate its cumulants to those of the sine process as $N\to \infty$. 
In particular, lemma~\ref{thm:kernel} which is  the main ingredient to prove proposition~\ref{thm:cumulant} below  is a straightforward adaptation of lemma~2.6~in~\cite{JL_15}.

\begin{lemma} \label{thm:kernel}
 We consider two families of kernels $(S_N)_{N>0}$ and $(\tilde S_N)_{N>0}$ on $\R$. If there exist $\beta>0$, $L>0$, and a function $\Gamma_N:\R\to\R^+$ such that when $N$ is sufficiently large:
\begin{itemize}
\item[(1)]    for all $x,y \in[-L,L]$,  $\displaystyle \big|\tilde S_N(x,y) -S_N(x,y) \big|  \le C_L N^{-\beta} $ . \\
\item[(2)]    for all $x,y \in[-L,L]$,   $\displaystyle \big| S_N(x,y) \big|  \le \Gamma_N(x-y) $ .\\
\item[(3)]   $\displaystyle  \int_{-2L}^{2L} \Gamma_N(s) ds \le C \log(L N) $ . 
\end{itemize} 
Then, for  all  $\epsilon>0$, $\ell\in\N$, and for  any functions $f_{N,1},\dots, f_{N,\ell}$ with support in $[-L,L]$ such that $\sup\big\{ \| f_{N,k} \|_\infty : k=1,\dots, \ell \} \le C_\ell$, we have      
\begin{equation*}  \tr[f_{N,1}\tilde S_N \cdots f_{N,\ell} \tilde S_N] =   \tr[f_{N,1}S_N \cdots f_{N,\ell} S_N]  +  \underset{N\to\infty}{O}\left( N^{-\beta+\epsilon}\right) \ .  \end{equation*}
\end{lemma}

We let 
\begin{equation*}  
 S_N(\xi,\zeta) =  \frac{\sin \left[ \pi N \big( F (x_0 + \xi N^{-\alpha} ) - F(x_0 + \zeta N^{-\alpha} )\big) \right]}{\pi(\xi- \zeta)} 
  \end{equation*}
and
\begin{equation} \label{L1B}
\Gamma_N(\xi-\zeta) = 
\begin{cases}  C _0N^{1-\alpha} &\text{if } |\xi-\zeta|^{-1}\le\frac{1}{N} \\
1/\pi |\xi-\zeta|  &\text{if } |\xi-\zeta|^{-1}>\frac{1}{N}
  \end{cases} \ .
\end{equation}
The asymptotics (\ref{sine}) implies that the kernels $\tilde K_N$ and $S_N$ satisfy condition (1) of  lemma~\ref{thm:kernel}. We claim that, if $C_0>0$ is sufficiently large, conditions (2) and (3) hold as well, so that we obtain for any $m_1,\dots, m_\ell \in \N$,
\begin{equation}  \label{cumulant_2}
\tr\left[f^{m_1}\tilde K_N \cdots f^{m_\ell} \tilde K_N \right] 
= \tr\left[f^{m_1} S_N \cdots f^{m_\ell} S_N \right] 
+  \underset{N\to\infty}{O}\left( N^{-\beta+\epsilon}\right) \ .
\end{equation}  
By (\ref{L1B}), it is straightforward  to check that for any $C_0>0$,  
\begin{equation*} 
\int_{-2L}^{2L} \Gamma_N(s) ds \le \log(L N)  + O(1) \ ,
\end{equation*}
so that condition (3) holds. To check condition (2), note that by definition of  $J$, (\ref{bulk}), for any $0<\epsilon_0<1/2$, there exists $\delta_0>0$ so that the density $\rho$ is continuous on  $[x_0-\delta_0,x_0+\delta_0]$ and for all $|x-x_0|<\delta_0$,
\begin{equation} \label{F'} 1- \epsilon_0 \le \
\frac{\rho(x)}{\rho(x_0)} \le 1+\epsilon_0  \ .
\end{equation}
If $N^\alpha> L/\delta_0 $ and  $C_0 \ge \rho(x_0)(1+\epsilon_0) $, this implies that for all $\xi,\zeta \in [-L,L]$,
$$ \left| F (x_0 + \xi N^{-\alpha} ) - F(x_0 + \zeta N^{-\alpha} )  \right|= N^{-\alpha} \left| \int_\xi^\zeta \rho(x_0 + s N^{-\alpha}) ds \right| \le C_0  N^{-\alpha} |\xi-\zeta|  \ .$$
Thus,  if we use the trivial bound $|\sin x| \le |x|\vee 1$, according to (\ref{L1B}), we conclude that 
$$ \big| S_N(\xi,\zeta) \big| \le \Gamma_N(\xi-\zeta) \ . $$


The map  $F$ is continuous non-decreasing, so it has a generalized inverse 
\begin{equation} \label{G}
G(x)= \inf \big\{ t \in \R : F(t) \ge x \big\} \ .
\end{equation}
In the sequel, we will assume that  $\delta_0$ is sufficiently small, so that (\ref{F'}) holds and the map  $G$ is continuously differentiable on $[F(x_0)-\delta_0, F(x_0)+\delta_0 ]$ with
\begin{equation} \label{G'}
G'(x)= \frac{1}{F'(G(x))} = \frac{1}{\rho(G(x))} \ .
\end{equation} 

Recall that the sine process $ \Xi^{\sin}_\nu$ is the determinantal process on $\R$ with a correlation kernel $ K^{\sin}_\nu$ given by (\ref{sine_K}).

\begin{proposition}  \label{thm:cumulant}
Let $f\in C_0(\R)$ and $ \alpha \in (0,1]$. We have for any $n\ge 1$, 
\begin{equation*}
\lim_{N\to \infty} \Cu^n\left[\Xi_N^{\alpha,x_0} f \right] 
=  \lim_{N\to \infty} \Cu^n\left[ \Xi^{\sin}_{\nu_N} f_N\right] \ ,
\end{equation*}
where  $\nu_N = N^{1-\alpha}\rho(x_0)$,
\begin{equation} \label{f_N}
f_N(x)= f\left( N^{\alpha} \left\{G\left( F\left(x_0\right)+  \rho(x_0) \frac{x }{N^{\alpha}} \right) -x_0 \right\}\right) \ ,
\end{equation}
 and the functions $F= F_\rho$  and $G$  is given by $(\ref{G})$.
\end{proposition}

\begin{remark}  \label{rk:f_N}
Observe that for any $0<\epsilon_0<1/2$, by $(\ref{F'})$ and $(\ref{G'})$, if $N^\alpha>2\rho(x_0) L /\delta_0$, then for all $x\in[-2L,2L]$,
\begin{equation*}
\frac{1}{1+\epsilon_0}\le \rho(x_0)G'\left( F\left(x_0\right)+  \rho(x_0) \frac{x }{N^{\alpha}} \right) \le \frac{1}{1-\epsilon_0} \ .
 \end{equation*}
If we integrate this estimate, since $f \in C_0\big([-L,L]\big)$, this implies that  the function $f_N \in C_0(\R)$ with support in $\big[-L_0, L_0\big]$ where $L_0=L(1+\epsilon_0)$. \\
\end{remark}

\noindent{\it Proof of proposition~\ref{thm:cumulant}.}
We  fix $m_1,\dots, m_\ell \in \N$ and we suppose that $N^\alpha>L (1 \vee 2\rho(x_0) ) /\delta_0$. We can make the change of variables 
\begin{equation} \label{unfold}
y_k = \frac{N^{\alpha}}{ \rho(x_0) }\left\{ F\left(x_0+ N^{-\alpha}x_k\right) - F\left(x_0\right) \right\}  
\end{equation}
in the formula
\begin{equation*}
 \tr\left[f^{m_1} S_N\cdots f^{m_\ell} S_N \right] 
=\int\limits_{-L}^L\cdots\int\limits_{-L}^L \prod_{k=1}^\ell f(x_k)^{m_k}\frac{\sin \left[ \pi N \big( ( F (x_0 + x_k N^{-\alpha} ) - F(x_0 + x_{k+1}N^{-\alpha} )\big) \right]}{\pi(x_k- x_{k+1})}
 d^\ell x \ .
\end{equation*}
If we let 
\begin{equation} \label{g_N}
g(y) = G\left(  F\left(x_0\right)+  \rho(x_0)  y   \right) -x_0  \ ,
\end{equation}
and $f_N$ be given by (\ref{f_N}), according to  remark~\ref{rk:f_N},  this leads to
\begin{align} \notag
 \tr\left[f^{m_1} S_N \cdots f^{m_\ell} S_N \right] 
&=\int\limits_{-L_0}^{L_0}\cdots\int\limits_{-L_0}^{L_0} \prod_{k=1}^\ell 
f_N(y_k)^{m_k} \frac{ g'( y_k N^{-\alpha}) \sin\left[ \pi \rho(x_0) N^{1-\alpha} \left( y_k-y_{k+1} \right)\right] }{\pi N^{\alpha}\big( g( y_k N^{-\alpha})-g(y_{k+1 }N^{-\alpha}) \big) } 
d^\ell y \\
& \label{cumulant_3}
=  \tr\left[f_N^{m_1} \tilde S_N \cdots f_N^{m_\ell} \tilde S_N \right]     \ , 
\end{align}
where
\begin{equation*}
\tilde S_N(x,y) : = \frac{ g'( y N^{-\alpha})\sin\left[ \pi \rho(x_0) N^{1-\alpha} \left( y-z \right)\right]}{N^{\alpha}\big\{ g( y N^{-\alpha})-g(zN^{-\alpha}) \big\}} \ .
\end{equation*}
For any $0< \alpha \le1$, a Taylor expansion in (\ref{g_N}) yields for all $y,z\in [-L_0,L_0]$,
\begin{equation*}  
 g'\left(yN^{-\alpha}\right)^{-1}N^{\alpha}\big\{ g\left( y N^{-\alpha}\right)-g\left(zN^{-\alpha}\right)\big\}
= (y-z) \left\{ 1+O\left( |y-z| N^{-\alpha}\right) \right\}  \ .
\end{equation*} 
This implies that uniformly for all $y,z\in [-L_0,L_0]$,
\begin{equation*} 
\tilde S_N(x,y) = \frac{\sin\left[ \pi \nu_N \left( y-z \right)\right]}{y-z} 
+O\left( N^{-\alpha}\right) \ ,
\end{equation*}
where  $\nu_N = N^{1-\alpha}\rho(x_0)$. Thus, the kernels $\tilde S_N$ and $K^{\sin}_{\nu_N}$  satisfy condition (1) of  lemma~\ref{thm:kernel} with $\beta=\alpha$. Moreover, if $\Gamma_N$ is given by (\ref{L1B}) with $C_0=\rho(x_0)$, the kernel   
$K^{\sin}_{\nu_N}$ also satisfies conditions (2) and (3). Therefore, since the functions $f_{N,k}=f_N^{m_k}  $ have support in $[-L_0,L_0]$ and 
$$\sup\big\{ \| f_{N,k} \|_\infty : k=1,\dots, \ell \} \le \|f\|_\infty^{m_1\vee\cdots\vee m_\ell} \ ,$$
 by lemma~\ref{thm:kernel}, we obtain
\begin{equation}  \label{cumulant_4}
\tr\left[f^{m_1}_N\tilde S_N \cdots f^{m_\ell}_N \tilde S_N \right] 
= \tr\left[f^{m_1}_N K^{\sin}_{\nu_N} \cdots f^{m_\ell}_N K^{\sin}_{\nu_N} \right] 
+  O\left( N^{-\alpha+\epsilon}\right) \ .
\end{equation}  
If we combine formulae (\ref{cumulant_2}), (\ref{cumulant_3}) and (\ref{cumulant_4}), we have proved that 
for any $m_1,\dots, m_\ell \in \N$,
\begin{equation}  \label{cumulant_5}
\tr\left[f^{m_1}\tilde K_N \cdots f^{m_\ell} \tilde K_N \right] 
= \tr\left[f^{m_1}_N K^{\sin}_{\nu_N} \cdots f^{m_\ell}_N K^{\sin}_{\nu_N} \right] 
+  O\left( N^{-\alpha\wedge\beta+\epsilon}\right) \ .
\end{equation}  
Since, by formula (\ref{cumulant_1}), the cumulants of the random variable $\Xi_N^{x_0,\alpha} f$ are linear combination of such traces, we conclude  by formula (\ref{cumulant_5}) that for any $n\ge 1$,
\begin{equation} \label{cumulant_6}
 \Cu^n\left[\Xi_N^{\alpha,x_0} f \right] 
=  \Cu^n\left[ \Xi^{\sin}_{\nu_N} f_N\right] +  O\left( N^{-\alpha \wedge \beta+\epsilon}\right) \ . 
\end{equation}\qed
 

\begin{remark} \label{rk:unfold}
In the physics literature, the change of variables $(\ref{unfold})$ is known as {\it unfolding the spectrum} since in the context of random matrices, it corresponds to rescaling the eigenvalue process  so that it has a constant density $\nu_N$ in a mesoscopic range around the point $x_0\in J_V$. Notice that in formula (\ref{sine}), if $\rho(x_0) \neq 0$, a Taylor expansion of the function $F_\rho$ shows that we recover the standard sine-kernel asymptotics in the regime $\alpha>1/2$,
 \begin{equation*} \label{sine_2}
\frac{1}{N^\alpha\varrho_V (x_0)} K_N\left(x_0 + \frac{\xi}{N^\alpha\varrho_V (x_0)} , x_0 + \frac{\zeta}{N^\alpha\varrho_V (x_0)} \right)= \frac{\sin \big[ \pi N^{1-\alpha} (\xi-\zeta)\big]}{\pi(\xi-\zeta)} + \underset{N\to\infty}{O}(N^{1-2\alpha}) \ .
\end{equation*}
Hence, at sufficiently small scales, the fact that the eigenvalues are not uniformly distributed is not relevant and, if the integrated density of states $F_\rho$ is smooth in $J_V$, we can deduce proposition~\ref{thm:cumulant} directly from lemma~\ref{thm:kernel} without  the change of variables $(\ref{unfold})$.\end{remark}

  
First, we use proposition~\ref{thm:cumulant} to derive the local correlations. In the regime $\alpha=1$,  for any $n\ge 1$,
\begin{equation} \label{cumulant_micro}
 \lim_{N\to\infty} \Cu^n\left[ \Xi_N^{1,x_0} f \right] = \lim_{N\to\infty} \Cu^n\left[ \Xi^{\sin}_{\rho(x_0)} f_N \right] \ . 
 \end{equation}
 By  (\ref{G'}), a Taylor expansion of the map $G$ yields for all $|x|<L_0$,
\begin{equation} \label{G_1}
  \lim_{N\to\infty} N^{\alpha} \left\{  G\left( F\left(x_0\right)+  \rho(x_0) \frac{x }{N^{\alpha}} \right) -x_0 \right\}= x \ .
  \end{equation}
By remark~\ref{rk:f_N}, the function $f_N$ has support in $[-L_0,L_0]$ and by continuity of $f$, the limit 
(\ref{G_1}) implies that $\displaystyle \lim_{N\to\infty} f_N(x)= f(x)$ for all $x\in\R$. Hence, by the dominated convergence theorem, we get
$$ \lim_{N\to\infty} \Cu^n\left[ \Xi^{\sin}_{\rho(x_0)} f_N \right]  =  \Cu^n\left[ \Xi^{\sin}_{\rho(x_0)} f \right] \ . $$
 By (\ref{cumulant_micro}), it proves that $\displaystyle  \lim_{N\to\infty} \Cu^n\left[ \Xi_N^{1,x_0} f \right]  = \Cu^n\left[ \Xi^{\sin}_{\rho(x_0)} f \right]$ for any $f\in C_0(\R)$ and the limit theorem (\ref{slt}) follows from the fact that the sine process is characterized by its cumulants. \\


We now turn to the proof of (\ref{clt}) in the mesoscopic regime, $0<\alpha<1$. The argument is different because, in formula (\ref{cumulant_6}), the density of the sine-process $\nu_N \to \infty$ as $N\to\infty$.
A relevant result in this regime is a CLT due to Soshnikov for the sine process.
\begin{theorem}[Thm.~4, \cite{Soshnikov_01}] \label{thm:Soshnikov} 
For any function $f\in H^{1/2}(\R)$, as $\nu\to\infty$,
$$\Xi^{\sin}_{\nu} f - \E{\Xi^{\sin}_\nu f}  \ \Rightarrow\ \No\left(0, \|f\|_{H^{1/2}}\right) \ . $$
\end{theorem}

The proof is based on Fourier analysis and a combinatorial argument given in the article~\cite{Soshnikov_00a}. Although the original proof is given for Schwartz functions, using a density argument, it is not  difficult  to extend Soshnikov's CLT to all test functions in the Sobolev space $H^{1/2}(\R)$. In order to deduce theorem~\ref{thm:universality} from proposition~\ref{thm:cumulant}, we see that it suffices to extend the proof of  
theorem~\ref{thm:Soshnikov} to deal with test functions $f_N$ of the form (\ref{f_N}).
 To proceed we need  further notations and to recall two key  lemmas from~\cite{Soshnikov_00a}.\\

For any tuple $m\in\N^\ell$, we define
\begin{equation} \label{MCL_function}
 \Upsilon_n(u_1,\dots,u_n)= \sum_{\ell=1}^n \frac{(-1)^{\ell+1}}{\ell} \sum_{\begin{subarray}{c} m_1, \dots, m_\ell \ge 1 \\ m_1+ \cdots + m_\ell = n  \end{subarray}  } \frac{n!}{m_1!\cdots m_\ell!}\  \underset{1\le i \le\ell}{\max}\{ u_1+\cdots+ u_{m_1+\cdots+m_i}\} 
 \end{equation}

\begin{lemma}[\cite{Soshnikov_00a}] \label{thm:S_lemma}  
There exists a constant $C_n>0$ which only depends on $n \ge 2$ such that  for any $\nu>0$ and any function $f\in L^1(\R)$,
\begin{equation*} \label{S_lemma}
\left|\Cu^n\big[\Xi^{sin}_\nu f\big] +2 
\int_{\R^n_0} \Re\bigg\{ \prod_i \hat{f}(u_i)  \bigg\} \Upsilon_n(u_1,\dots,u_n) d^{n-1}u \right| \le C_n \int_{\mathcal{A}^n_\nu} \bigg|\prod_i \hat{f}(u_i)  \bigg| \big(|u_1|+\cdots |u_n|\big) d^{n-1} u \ ,   
\end{equation*}
where $\R^n_0=\big\{u\in\R^n : u_1+\cdots+u_n=0\big\}$ and $\mathcal{A}^n_\nu= \big\{ u\in\R^n_0 :  |u_1|+\cdots +|u_n| > \nu \big\}$.
\end{lemma}

\begin{lemma}[Main Combinatorial lemma, \cite{Soshnikov_00a}] \label{thm:MCL}
For any $u\in\R^n_0$,
\begin{equation*} \label{MCL}
 \sum_{\sigma\in \Sy_n} \Upsilon_n\left(u_{\sigma(1)},\dots,u_{\sigma(n)}\right) 
= \begin{cases} |u_1|  &\text{if } n=2 \\
 0  &\text{if } n>2  \end{cases} \ .
 \end{equation*}
\end{lemma}


If $g\in C^1_0(\R)$, we define 
\begin{equation} \label{norm_3}
\| g\|_{H^1}^2 = \int_\R |\hat g (u)|^2 |u|^2 du = \frac{1}{4\pi^2}  \int_\R \big| g'(x) \big|^2 dx \ .
\end{equation}
We will also need the following result.

\begin{lemma}  \label{thm:f_N'}
If $f\in C^1_0(\R)$ and the function $f_N$ is given by $(\ref{f_N})$, then
\begin{equation*}
 \lim_{N\to\infty} \left\| f_N-f\right\|_{H^{1}} = 0 \ .
\end{equation*}
\end{lemma}

\proof Since $G\in C^1\big([F(x_0)-\delta_0, F(x_0)+\delta_0 ]\big)$, by   remark~\ref{rk:f_N},   if  $N^\alpha>2\rho(x_0) L /\delta_0$,
the functions $f_N$ are continuously differentiable on $\R$ and 
\begin{equation*} 
f_N'(x)= \rho(x_0) G'\left( F\left(x_0\right)+  \rho(x_0) \frac{x }{N^{\alpha}} \right)  f'\left( N^{\alpha} \left\{G\left( F\left(x_0\right)+  \rho(x_0) \frac{x }{N^{\alpha}} \right) -x_0 \right\}\right) \ .
\end{equation*}
Then, by the triangle inequality,
\begin{align} 
\label{f_N'}
\big| f_N'(x)-f'(x) \big| &\le  \|f'\|_\infty  \left| \rho(x_0) G'\left( F\left(x_0\right)+  \rho(x_0) \frac{x }{N^{\alpha}} \right) -1 \right| \\
&\notag\hspace{1cm}
+ \left|f'\left( N^{\alpha} \left\{G\left( F\left(x_0\right)+  \rho(x_0) \frac{x }{N^{\alpha}} \right) -x_0 \right\}\right)- f'(x) \right| \ .
\end{align}
First note that, by (\ref{G_1}) and the continuity of $f'$,
\begin{equation}\label{limit_2}
 \lim_{N\to\infty}   \left|f'\left( N^{\alpha} \left\{G\left( F\left(x_0\right)+  \rho(x_0) \frac{x }{N^{\alpha}} \right) -x_0 \right\}\right)- f'(x) \right| =0 \ .
\end{equation}
Second, by  remark~\ref{rk:f_N},  we have for all sufficiently large $N$,
\begin{equation*}
\left|\rho(x_0)G'\left( F\left(x_0\right)+  \rho(x_0) \frac{x }{N^{\alpha}} \right) -1  \right| \le \frac{\epsilon_0}{1-\epsilon_0} \ .
 \end{equation*}
Since $\epsilon_0$ can be taken arbitrary small, we deduce that
\begin{equation}\label{limit_3}
 \lim_{N\to\infty}\left| \rho(x_0) G'\left( F\left(x_0\right)+  \rho(x_0) \frac{x }{N^{\alpha}} \right) -1 \right| 
=0 \ . 
\end{equation}
In the end, by the dominate convergence and the estimates  (\ref{f_N'} - \ref{limit_3}), we conclude that
$$  \lim_{N\to\infty} \left\| f_N-f\right\|_{H^{1}}  = \lim_{N\to\infty} \frac{1}{4\pi^2} \int_{-L_0}^{L_0} \big| f_N'(x)-f'(x) \big|  = 0 \ .$$ \qed

Observe that, if $g\in C_0^1(\R)$, according to (\ref{norm_1}) and (\ref{norm_3}),
\begin{align*} \| g \|_{H^{1/2}}^2 &\le \| \hat g \|_\infty^2 +  \| g \|_{H^1}^2 \\
& \le  \| g \|_{L^1}^2 +  \| g \|_{H^1}^2  \ .
\end{align*}
By (\ref{G_1}) and the dominated convergence theorem, we get $\displaystyle \lim_{N\to\infty}  \| f_N-f \|_{L^1} = 0$. Thus,  by lemma~\ref{thm:f_N'}, we obtain for any $f\in C_0^1(\R)$,
\begin{equation} \label{limit_0}
 \lim_{N\to\infty} \left\| f_N-f\right\|_{H^{1/2}} = 0 \ .
\end{equation}


For now, let us also assume that, with  $\nu_N=N^{\alpha}\rho(x_0)$ and $\mathcal{A}^n_\nu$  defined in  lemma ~\ref{thm:S_lemma}, we have
\begin{equation} \label{limit_1}
\lim_{N\to\infty} \int_{\mathcal{A}^n_{\nu_N}} \bigg|\prod_i \hat{f}_N(u_i)  \bigg| \big(|u_1|+\cdots |u_n|\big) d^{n-1} u = 0 \ .
\end{equation}
This implies that  for any $n\ge 2$, 
\begin{equation*}
\lim_{N\to\infty} \Cu^n\left[ \Xi^{\sin}_{\nu_N} f_N\right] 
= 2 \lim_{N\to\infty} \int_{\R^n_0} \Re\bigg\{ \prod_i \hat{f}_N(u_i)  \bigg\} \Upsilon_n(u_1,\dots,u_n) d^{n-1}u \ .
\end{equation*}
Since $f$ is real-valued and $\Upsilon_2(u,-u)=|u|/2$,  by lemma ~\ref{thm:MCL}, we get
\begin{equation*}
\lim_{N\to\infty} \Cu^n\left[ \Xi^{\sin}_{\nu_N} f_N\right] 
= \begin{cases} \displaystyle
\lim_{N\to\infty}   \int_\R \big|\hat f_N(u) \big|^2|u|  du    &\text{if } n=2 \\
 0  &\text{if } n>2  \end{cases} \ .
 \end{equation*}
By  proposition~\ref{thm:cumulant}  and (\ref{limit_0}), we conclude that  for any $f\in C_0^1(\R)$,
 \begin{equation*}
\lim_{N\to \infty} \Cu^n\left[\Xi_N^{\alpha,x_0} f \right] =
 \begin{cases} \|f\|_{H^{1/2}}^2 &\text{if } n=2 \\
 0  &\text{if } n>2  \end{cases} \ .
\end{equation*}
 
 A special case of the limit (\ref{limit_1}) was computed in \cite[proposition~4.13]{JL_15}.
 The proof relies on lemma~\ref{thm:f_N'} and it is straightforward to generalize the argument of \cite{JL_15} to obtain (\ref{limit_1}).
Hence, by (\ref{cumulant_limit}), the CLT (\ref{clt}) holds for any $0<\alpha<1$, $x_0\in J$, and $f\in C_0^1(\R)$.

\section{Unitary invariant Ensembles} \label{sect:RMT}

\subsection{General context} \label{sect:intro}


The most well-known probability measure on the space of $N\times N$ Hermitian matrices is the Gaussian Unitary Ensemble,
\begin{equation}\label{GUE}
 d\mathbb{P}^{\G}_N = \Z_{N}^{-1}   e^{- N \tr H^2} dH \ , 
 \end{equation}
where $dH$ denotes the Lebesgue measure.  In this section, we will consider some generalizations of the GUE  of the form
\begin{equation}\label{UE}
 d\mathbb{P}^\omega_N =  \Z_{\omega,N}^{-1}e^{\tr  \log \omega(H)} dH \ ,
 \end{equation}
where the function $\omega:\R \to [0,+\infty)$ is upper-semicontinuous and such that for all $k\ge 0$, 
\begin{equation} \label{weight}   \int |x|^{k} \omega(x) dx < \infty \ . 
\end{equation}

This condition implies that the partition function $ \Z_{\omega,N}  <\infty$.
For scaling reasons, the weight $\omega$ may also depend on the dimension $N$ even though we will not  indicate it to keep our notations as simple as possible.
The matrix $\log \omega(H)$ is defined by functional calculus and the trace guarantees that the measure (\ref{UE}) is invariant under the transformation $H \mapsto UHU^*$ for any $U\in \mathcal{U}(N)$.
Hence, the name unitary invariant ensembles. In particular, if we use the spectral decomposition of $H$, under $\mathbb{P}^\omega_N$, the eigenvectors are independent  of the  spectrum $\Lambda$ and $\Lambda=\big\{\lambda_1,\dots, \lambda_N\big\}$ has a joint density on $\R^N$ which is given by  
\begin{equation}\label{density}
 \mathscr{F}^\omega_N(x_1,\dots, x_N) = Z_{\omega,N}^{-1} \det_{N\times N}\big[ x_j^{k-1} \big]   \det_{N\times N}\big[ x_j^{k-1} \omega(x_j)\big] \ .
 \end{equation}
 
 
In order to analyze the probability measure $\mathbb{P}^\omega_N$, a method introduced by  Gaudin and Metha in~\cite{MG_60} consists in rewriting the density $ \mathscr{F}^\omega_N$ using the orthogonal polynomials with respect to the measure $\omega(x) dx$ on $\R$. The condition (\ref{weight}) guarantees that these polynomials exist and we define
 for any $k\ge 0$,
\begin{equation} \label{OP}
 \pi_k(x) = x^k + \alpha_k x^{k-1}+ \cdots 
  \hspace{1cm}\text{and}\hspace{1cm} 
\int \pi_k(x) \pi_j(x) \omega(x) dx  = \gamma_k^{-2} \1_{k=j}    \ .
\end{equation}

Then, it follows from formula (\ref{density}) that  the eigenvalues density is
\begin{equation}\label{density_OP}
 \mathscr{F}^{\omega}_N(x_1,\dots, x_N) = \frac{1}{N!} \det_{N\times N}\big[ K_N^\omega(x_j, x_k)  \big] dx_1\cdots dx_N    \ ,
 \end{equation}
where 
\begin{equation}\label{kernel_CD}
K_N^\omega(x, y) = \gamma_{N-1}^2\frac{ \pi_N(x) \pi_{N-1}(y)- \pi_{N-1}(x) \pi_{N}(y)}{x-y}\sqrt{\omega(x)\omega(y)}   \ .
\end{equation}

Formulae (\ref{density_OP} - \ref{kernel_CD}) implies that  the point process $\Lambda$ is  determinantal with correlation kernel  $K_N^\omega$ in the sense of (\ref{det}). 
These facts are well-established and we refer to e.g$.$ \cite{Deift_99, Konig_05} for an introduction to the subject. 
%
%
%
 %
By theorem~\ref{thm:universality},  this reduces the question of universality of mesoscopic fluctuations for the ensembles (\ref{UE})  to obtain a precise asymptotics for the OPs with respect to the measure   $\omega(x) dx$. \\

 Beyond the  context of random matrix theory, one may consider the determinantal process (\ref{density_OP}) associated with a general measure. These processes are  known as orthogonal polynomial ensembles and a significant amount of research has focused on proving the sine-process universality at the local scale, see \cite{Lubinsky_09b} and reference therein. At mesoscopic scales, another universality result just  appeared in~\cite{BD_15} and the authors already obtained theorem~\ref{thm:JE} below.
 Instead of working with the correlation kernel of the process, they reformulate the cumulant problem in terms of the so-called Jacobi matrix of the measure $\omega(x)dx$ and this reduces the question of universality to the asymptotics of the recurrence coefficients which define the OPs. The drawback of their method is that, for technical reasons, it fails when the reference measure depends on the dimension $N$, like the GUE or the exponential weights considered in section~\ref{sect:UIE}. However, this other method requires only the asymptotics of the recurrence coefficients  and 
 it applies to discrete or singular measures where the OP asymptotics is difficult to derive. \\

Under general conditions and provided that the weight $\omega$ is suitably normalized as $N\to\infty$, see \cite{BD_14,Hardy_15}, it is known that  there is a Law of Large numbers
\begin{equation} \label{equilibrium'}
 \frac{1}{N} \sum_{k=1}^N f(\lambda_k) \underset{N\to\infty}{\longrightarrow}  \int f(x) d\mu_\omega(x) dx    \ , \hspace{1cm}  \mathbb{P}^\omega_N \text{ - almost surely.}
\end{equation}
Hence, $\mu_\omega$ is the equilibrium measure for the ensemble $\mathbb{P}^\omega_N$ and it  has compact support. 
In the following, we will suppose that it is absolutely continuous:  $d\mu_\omega=\varrho_\omega(x) dx$.
 The equilibrium density $\varrho_\omega$ plays a fundamental role in the non-linear steepest descent  introduced in \cite{DZ_93} and it comes in the asymptotics of the Christoffel-Darboux  kernel.  Based on the results of \cite{Deift_al_99_b,  KV_02, Kuijlaars_al_04}, in sections~\ref{sect:UIE} and~\ref{sect:JE}, we will derive the mesoscopic asymptotics for the correlation kernels of the ensembles $\mathbb{P}^V_N$, and the so-called  modified Jacobi ensembles respectively. 
 We do not intend to  review the Riemann-Hilbert literature but we point out that the Deift-Zhou steepest descent   has been developed by several authors and it has yields local universality for an extensive pool of OP ensembles and it should be possible to apply theorem~\ref{thm:universality} to prove mesoscopic universality as well, e.g.~for the modified Laguerre ensembles and Wishart matrices using the results of~\cite{Vanlessen_07}. \\

In section~\ref{sect:variance}, assuming that the OPs  satisfy {\it classical} asymptotic formulae, see (\ref{semiclassical_1}) below, we derive estimates for the variance of linear statistics. This allows us to extend the CLT (\ref{clt}) to a larger class of test functions, 
 see theorems~\ref{thm:clt_1} and~\ref{thm:JE} below. Lemma~\ref{thm:variance} is also of interest  for global linear statistics $(\alpha=0)$. In particular, in the appendix~\ref{A:variance}, we extend the scope of the CLT (\ref{clt_2}) to rather general test functions.
Finally, for the GUE, it is possible to get the mesoscopic asymptotics of the correlation kernel without using the Riemann-Hilbert techniques; a complete proof is given in  section ~\ref{sect:GUE}.

\subsection{Variance estimates} \label{sect:variance}

For the GUE, estimates for the variance of mesoscopic linear statistics have been derived in~\cite{FKS_13} for Schwartz-class test functions.  Using a similar formalism, we will derive estimates for the variance of linear statistics which are valid both  in the global and mesoscopic regimes for arbitrary weight $\omega$  such that the support of the equilibrium density $\varrho_\omega$ is  connected. Our method relies only on  the bulk asymptotics of the OPs, therefore it cannot yield optimal upper-bounds. However, our results  apply to test functions with rather mild smoothness and slow  decay, such as the functions $g_t$ which arise when considering the logarithm of a regularized characteristic polynomial, cf.~(\ref{W0}) below. \\

  According to formula (\ref{OP}), we let $\Phi_k(x)= \gamma_k \pi_k(x) \sqrt{\omega(x)}$ for any $k\ge 0$. Thus $\big(\Phi_k\big)_{k=0}^\infty$ is an orthonormal family in $L^2(\R)$.\
 By (\ref{kernel_CD}), the Christoffel-Darboux kernel for the weight $\omega(x)$ on $\R$  is given by
\begin{equation} \label{kernel_P}
K_N^\omega(x, y) = \frac{\gamma_{N-1}}{\gamma_{N}}\frac{ \Phi_N(x) \Phi_{N-1}(y)- \Phi_{N-1}(x) \Phi_{N}(y)}{x-y}   \ .
\end{equation}
We suppose that the OPs have the following asymptotics for all  $|x|<1$,
\begin{equation} \label{semiclassical_1}
 \Phi_N(x)= \sqrt{\frac{2}{\pi}} \frac{\cos\big[ N\pi F(x) + \psi(x)\big]}{(1-x^2)^{1/4}} + \underset{N\to\infty}{o_\epsilon(1)} \ .
\end{equation}
The function $\psi \in C(-1,1)$ and the notation $o_\epsilon$ means that the error term converges to 0 uniformly  for all $|x|<1-\epsilon$, i.e.~it only depends on the parameter $\epsilon>0$. Moreover, according to (\ref{F}), $F=F_{\varrho_\omega}$  where $\varrho_\omega$ denotes the equilibrium density for the ensemble $\mathbb{P}^\omega_N$ and we assume that $J_{\varrho_\omega} = (-1,1)$.  
For a fixed weight $\omega$, we can deduce from formula (\ref{semiclassical_1}) that
\begin{equation} \label{semiclassical_2}
 \Phi_{N-1}(x)= \sqrt{\frac{2}{\pi}} \frac{\cos\big[ N\pi F(x) + \tilde\psi(x)\big] }{(1-x^2)^{1/4}} 
 + \underset{N\to\infty}{o_\epsilon(1)} \ ,
\end{equation}
where $ \tilde\psi(x) = \psi(x) -\pi F(x)$.  When the weight $\omega$ depends on the dimension $N$, we will suppose that formula (\ref{semiclassical_2})  still hold for some function $\tilde\psi \in C(-1,1)$.
For instance, when $\omega(x)=e^{-N V(x)}$ and $V(x)$ is analytic on $\R$, the asymptotics of the OPs   have been investigated in~\cite{Deift_al_99_b} by solving the appropriate Riemann-Hilbert problem and it is rather straightforward to check that, if $J_V=(-1,1)$, then the asymptotics (\ref{semiclassical_1} - \ref{semiclassical_2}) hold, $\psi(x)= -\tilde\psi(x) = \arcsin(x)/2$ regardless of the potential $V(x)$, and
furthermore
\begin{equation}\label{semiclassical_3}
 \lim_{N\to\infty}  \frac{\gamma_{N-1}}{\gamma_{N}} = \frac{1}{2} \ .
\end{equation}
For the GUE, these results follows directly from the Plancherel-Rotach asymptotics, see section~\ref{sect:GUE}. \\

We consider the determinantal process $\Xi_N$ with correlation kernel (\ref{kernel_P}) and for any continuous function, we denote  $\Xi_N f  := \sum f (\lambda_k)$ where the sum is over the configuration $\{\lambda_k\}_{k=1}^N$. It is well-known that since  $K_N^\omega$ defines a projection on $L^2(\R)$, cf.~e.g.~\cite[Lem.~3.1]{JL_15}, we have 
\begin{equation}\label{variance_1}
 \Var\left[ \Xi_N f \right] = \frac{1}{2} \iint \left|f(x)-f(y) \right|^2  \big|K^\omega_N(x,y)\big|^2 dxdy \ .
\end{equation}

We have seen in the introduction that, for the GUE  or a general ensemble $\mathbb{P}^V_N$  satisfying the condition (\ref{1cut}), the CLT (\ref{clt_2}) implies that for any sufficiently smooth function $f$, $\Var\left[ \Xi_N f \right] \rightarrow \Sigma(f)^2$ as $N\to\infty$.  The question, we address in this section is wether there exists a constant $\mathrm{C}>0$ which may depend only on the weight $\omega$ such that  for more general test functions,
\begin{equation} \label{variance_0}
 \Var\left[ \Xi_N f \right] \le \mathrm{C}\  \Sigma(f)^2 \ . 
 \end{equation}

Since the point process $\Xi_N$ is essentially supported in the bulk $J$, we expect that, apart form some reasonable growth assumption,  the behavior of the function $f(x)$ outside $J$ should be irrelevant to  estimate $ \Var\left[ \Xi_N f \right]$. However, because of the effect of the spectral edge, it is technical to prove (\ref{variance_0}) in general. Instead, we will show that, if the OPs satisfy the conditions (\ref{semiclassical_1} - \ref{semiclassical_3}), then for any function $f\in H^{1/2}(\R)$ such that there exists $0<\delta<1$ and $\mathrm{L}>0$ so that
\begin{equation} \label{L} \tag{H.1}
 \sup \left\{ \left| \frac{f(x)- f(y)}{x-y} \right|  :  |x|\vee|y| > 1-\delta  \right\} < \mathrm{L} \ ,
 \end{equation}
we have
\begin{equation} \label{variance_global} 
 \limsup{N\to\infty}  \Var\left[ \Xi_N f \right]  \le  8 \iint\limits_{[-1,1]^2} \left| \frac{f(x)- f(y)}{x-y} \right|^2 \frac{dxdy}{\sqrt{1-x^2}\sqrt{1-y^2}}  \ . 
 \end{equation}

Note that, if $f\in H^{1/2}(\R)$, the condition (\ref{L}) guarantees that 
\begin{equation} \label{variance_2} 
\tilde{\Sigma}(f)^2 := \frac{1}{\pi^2}  \iint\limits_{[-1,1]^2} \left| \frac{f(x)- f(y)}{x-y} \right|^2 \frac{dxdy}{\sqrt{1-x^2}\sqrt{1-y^2}} <\infty  \ ,
 \end{equation} 
and $\Sigma(f) \le  \tilde\Sigma(f)$. 
In fact, if it exists, it is difficult to exhibit a function $ h\in H^{1/2}(\R)$ such that $\Sigma(h)<\infty$ and  $\tilde\Sigma(h)=\infty$. We begin by proving a simple lemma on the asymptotic behavior of the $L^2$-mass of the function $\Phi_N$.

\begin{lemma} \label{thm:localization} Suppose that formula $(\ref{semiclassical_1})$ holds, then for any $0<\epsilon<1$,  if $J_\epsilon= [-1+\epsilon, 1-\epsilon]$, we have 
\begin{equation*}
\int_{J_\epsilon} \big|\Phi_{N}(x)\big|^2 dx=\frac{2\arcsin(1-\epsilon)}{\pi}  + \underset{N\to\infty}{o_\epsilon(1)}  \ . 
\end{equation*}
\end{lemma}

\proof By (\ref{semiclassical_1}), for any $|x|<1$,
$$\big|\Phi_{N}(x)\big|^2=\frac{2}{\pi\sqrt{1-x^2}}\bigg\{ \frac{1+ \cos \big[2N\pi F(x) + 2\psi(x)\big]}{2}  
+ o_\epsilon(1)\bigg\} \ , $$
and this implies that 
$$\int_{J_\epsilon}\big|\Phi_{N}(x)\big|^2 dx= \int_{J_\epsilon} \frac{1+ o_\epsilon(1)}{\pi\sqrt{1-x^2}}dx + \int_{J_\epsilon} \frac{\cos\big[2N\pi F(x) + 2\psi(x)\big]}{\pi\sqrt{1-x^2}}dx \ . $$
The first integral converges to $2\arcsin(1-\epsilon)/\pi$ as $N\to\infty$ and it remains to show that the second integral which is oscillatory converges to 0.
By assumption, $F'(x)=\varrho_\omega(x)>0$ for all $x\in (-1,1)$ and we can make the change of variable $x= G(y)$ where $G=F^{-1}$, (\ref{G}). So that, if $g\in L^1(\R)$ with compact support in $(-1,1)$, we have  
$$ \int g(x) e^{i2N\pi F(x)} dx = \int g(G(y)) G'(y) e^{i2N\pi y} dy \ ,    $$
and, since $y\mapsto g(G(y)) G'(y)$ is integrable,
$\displaystyle\lim_{N\to\infty}  \int g(x) e^{i2N\pi F(x)} dx =0$ by the Riemann-lebesgue lemma.  
 Applying this argument to the function $g(x)=\frac{e^{i2\psi(x)}}{\pi\sqrt{1-x^2}} \1_{J_\epsilon}(x)$, we conclude that for any $0<\epsilon<1$,
$$\lim_{N\to\infty}\int_{J_\epsilon} \frac{\cos\big[2N\pi F(x) + 2\psi(x)\big]}{\pi\sqrt{1-x^2}}dx = 0  $$
and the proof is complete. \qed\\

According to formula (\ref{variance_1}), in order to get an upper-bound for the variance of linear statistics, we need to estimate  the quantity $|K^\omega_N(x,y)|^2$ for all $x,y \in \R$.  By (\ref{kernel_P}), we have
\begin{equation}  \label{variance_3} 
\big|K^\omega_N(x,y)\big|^2 \le  \frac{2\gamma_{N-1}}{\gamma_{N}}\frac{\big|\Phi_N(x) \Phi_{N-1}(y)\big|^2+\big| \Phi_{N-1}(x) \Phi_{N}(y)\big|^2}{|x-y|^2} \ ,
\end{equation}
Moreover, using the asymptotics (\ref{semiclassical_1} - \ref{semiclassical_3}), we get for all $|x|,|y|<1-\epsilon$,   
\begin{equation}  \label{variance_4} 
\big|K^\omega_N(x,y)\big|^2 \le \frac{1}{\pi^2} \frac{ \Delta_N(\epsilon)}{\sqrt{1-x^2}\sqrt{1-y^2}} \frac{1}{|x-y|^2} \ ,
\end{equation}
where $\Delta_N(\epsilon)>0$ and $\displaystyle\lim_{N\to\infty} \Delta_N(\epsilon) = 8$ for any $0<\epsilon<1$.
Using these two estimates and lemma~\ref{thm:localization} we can prove formula (\ref{variance_global}).
Lemma~\ref{thm:variance} is formulated for test functions in the Sobolev space  $H^{1/2}(\R)$. The smoothness condition is plainly necessary. However, by exploiting the decay of the Christoffel-Darboux kernel outside of the bulk, the result holds for test functions with polynomial growth as well, see appendix~\ref{A:variance}.

\begin{lemma} \label{thm:variance}
Suppose that the OPs with respect to the weight $\omega(x)$ on $\R$ satisfy the conditions $(\ref{semiclassical_1} - \ref{semiclassical_3})$, then for any function $f\in H^{1/2}(\R)$ which satisfies $(\mathrm{\ref{L}})$, we have for any $0<\epsilon<\delta$,
\begin{equation*}
  \Var\left[ \Xi_N f \right]  \le  \Delta_N(\epsilon)\tilde\Sigma(f)^2 + 
  \underset{N\to\infty}{O}\big( \mathrm{L}^2 \Theta(\epsilon)\big)  \ ,
 \end{equation*}
 where $\Theta(\epsilon)=1-\frac{2\arcsin(1-\epsilon)}{\pi}$,  $\Delta_N(\epsilon)>0$ and $\displaystyle\lim_{N\to\infty} \Delta_N(\epsilon) = 8$ for any $0<\epsilon<1$.
\end{lemma}

\proof  Fix $0<\epsilon<1$ and let $J_\epsilon= [-1+\epsilon, 1-\epsilon]$. By formulae (\ref{variance_2}) and (\ref{variance_4}), we have
\begin{equation}  \label{variance_5}
 \iint_{J_\epsilon^2} \left|f(x)-f(y) \right|^2  \big|K_N^{\omega}(x,y)\big|^2 dxdy \le \Delta_N(\epsilon) \tilde\Sigma(f)^2 \ .
\end{equation}
On the other hand, if $f$ satisfies the hypothesis (\ref{L}) and $0<\epsilon<\delta$, by formula (\ref{variance_3}),
\begin{equation} \label{variance_6}
 \iint_{\R^2\backslash J_\epsilon^2} \left|f(x)-f(y) \right|^2  \big|K_N^{\omega}(x,y)\big|^2 dxdy 
 \le \frac{4\gamma_{N-1}}{\gamma_{N}} \mathrm{L}^2  \iint_{\R^2\backslash J_\epsilon^2} \big| \Phi_{N}(x) \Phi_{N-1}(y)\big|^2 dx dy \ \ .
\end{equation}
By symmetry
\begin{align*}
 \iint\limits_{\R^2\backslash J_\epsilon^2} \big| \Phi_{N}(x) \Phi_{N-1}(y)\big|^2 dx dy
\le \int_\R  \big|\Phi_{N-1}(y)\big|^2 dy \int_{\R\backslash J_\epsilon} \big| \Phi_{N}(x)\big|^2 dx   
+  \int_\R  \big|\Phi_{N}(x)\big|^2 dx \int_{\R\backslash J_\epsilon} \big| \Phi_{N-1}(y)\big|^2 dy \ ,
\end{align*}
and since $\|\Phi_N\|_{L^2}= \|\Phi_{N-1}\|_{L^2} =1$, by lemma~\ref{thm:localization}, we obtain
\begin{equation*}
\limsup{N\to\infty}\iint_{\R^2\backslash J_\epsilon^2} \big| \Phi_{N}(x) \Phi_{N-1}(y)\big|^2 dx dy
\le 2\left( 1- \frac{2\arcsin(1-\epsilon)}{\pi} \right) \ . 
\end{equation*}
Note that we used that the asymptotics of lemma~\ref{thm:localization} holds for the function $\Phi_{N-1}$ as well; this follows from formula (\ref{semiclassical_2}).
Since $ \frac{\gamma_{N-1}}{\gamma_{N}} \to \frac{1}{2}$, this upper-bound and (\ref{variance_6}) implies that for any $0<\epsilon<\delta$,
\begin{equation} \label{variance_7}
 \iint_{\R^2\backslash J_\epsilon^2} \left|f(x)-f(y) \right|^2  \big|K_N^{\omega}(x,y)\big|^2 dxdy =
 \underset{N\to\infty}{O} \big(\mathrm{L}^2\Theta(\epsilon)\big) \ .
\end{equation}
The claim follows from formula (\ref{variance_1}) by combining the estimates (\ref{variance_5}) and (\ref{variance_7}). \qed\\

The first consequence of lemma~\ref{thm:variance} is that, for any function $f\in H^{1/2}(\R)$ which satisfies $(\mathrm{\ref{L}})$, there exists a constant $C>0$ such that  for any $0<\epsilon<1$,
\begin{equation*}
 \limsup{N\to\infty} \Var\left[ \Xi_N f \right]  \le 8\ \tilde\Sigma(f)^2 + 
C\left(1-\frac{2\arcsin(1-\epsilon)}{\pi}\right)  \ .
 \end{equation*}
Since the l.h.s.~is independent of $\epsilon$ and ${\displaystyle \lim_{\epsilon\searrow0} } \left( 1- \frac{2\arcsin(1-\epsilon)}{\pi} \right) =0$, this implies formula (\ref{variance_global}).
In the remainder of this section, we discuss the implication of lemma~\ref{thm:variance} for mesoscopic linear statistics.

\begin{proposition} \label{thm:variance_meso}
Suppose that the OPs with respect to the weight $\omega(x)$ on $\R$ satisfy the conditions $(\ref{semiclassical_1} - \ref{semiclassical_3})$. For any function $f\in H^{1/2}(\R)$ such  there exists $\mathrm{L}>0$ and
\begin{equation}  \label{L'} \tag{H.2}
\limsup{|x|\to\infty} \sup\left\{  |x|    \left| \frac{f(x)- f(y)}{x-y} \right| : |y| \le |x| \right\} < \mathrm{L} \ ,
\end{equation}
 for any $|x_0|<1$ and for any $0<\alpha<1$,  we have
\begin{equation} \label{variance_meso}
 \limsup{N\to\infty} \Var\left[ \Xi_N^{x_0,\alpha} f \right]  \le 32\ \| f\|_{H^{1/2}}^2 \ .
 \end{equation}
\end{proposition}

\begin{remark}
The main difficulty to estimate the variance of linear statistics is to control  the contribution from the edges of the spectrum. An issue that we avoided by using lemma~\ref{thm:localization} and  the condition 
$(\operatorname{\ref{L'}})$. Based on the results of \cite{Deift_al_99_b}, the same method should apply to the ensemble $\mathbb{P}_N^V$ in the multi-cut case as well, though the asymptotics of the OPs is more complicated and the argument becomes rather technical. It is straightforward to check that $(\operatorname{\ref{L'}})$ holds in both cases:\\
\phantom{}\hspace{.5cm}  i) $f\in C^1(\R)$ and $|f'(x)| \le \mathrm{L}|x|^{-1}$.\\
\phantom{}\hspace{.5cm}  ii) $f$ is bounded and has compact support.\\
In particular, the estimate $(\ref{variance_meso})$ applies to the resolvent $x\mapsto (x-z)^{-1}$  for any $z\in \C$ such that $\Im z \neq 0$ and for any function in $H^{1/2}\cap L^\infty(\R)$ with compact support.
 From the point of view of mesoscopic linear statistics, this encompass the most relevant class of test functions. 
 \end{remark}

\noindent {\it Proof of proposition~\ref{thm:variance_meso}.}  
The assumption (\ref{L'}) implies that there exists $\mathrm{C}>0$ so that, if $|x|\ge \mathrm{C}$, then for all $|y| \le |x|$,
\begin{equation*} 
  \left| \frac{f(x)- f(y)}{x-y} \right| \le \frac{\mathrm{L}}{|x|} \ .
  \end{equation*}
If we let $ g_N(x)=f\big( N^{\alpha}(x-x_0)\big)$,  we obtain for all $|x-x_0|> \mathrm{C}N^{-\alpha} $ and all $|y-x_0|<|x-x_0|$,
$$
\left| \frac{g_N(x)- g_N(y)}{x-y} \right| \le \frac{\mathrm{L}}{|x-x_0|} \ .
$$
This inequality shows that, if $|x_0-1| \wedge |x_0+1| =2\delta$, then for all $N> ( \mathrm{C}/\delta)^{1/\alpha}$  and for all $|x|>1-\delta$,
\begin{equation} \label{variance_10}
\left| \frac{g_N(x)- g_N(y)}{x-y} \right| \le \frac{\mathrm{L}}{\delta} \ .
\end{equation}
Hence, since the r.h.s.~of (\ref{variance_10}) is symmetric in $x$ and $y$, for sufficiently large $N$, the functions $g_N$ satisfy the condition (\ref{L})  and by lemma~\ref{thm:variance}, 
\begin{equation} \label{variance_11}
  \Var\left[ \Xi_N^{x_0,\alpha} f \right] = \Var\left[ \Xi_N g_N \right]   \le  \Delta_N(\epsilon)\tilde\Sigma(g_N)^2 + 
  \underset{N\to\infty}{O}\big(\Theta(\epsilon)\big)  \ .
 \end{equation}
By formula (\ref{variance_2}), if we let $J_\epsilon= [-1+\epsilon, 1-\epsilon]$ for $0<\epsilon<1$, then 
\begin{align} \label{variance_12}
&\tilde{\Sigma}(g_N)^2= \mathrm{I}_1(f ; N , \epsilon) +  \mathrm{I}_2(f ; N, \epsilon) \\
&\notag =  \frac{1}{\pi^2}  \iint\limits_{J_\epsilon} \left| \frac{g_N(x)- g_N(y)}{x-y} \right|^2  \frac{dxdy}{\sqrt{1-x^2}\sqrt{1-y^2}} 
+ \frac{1}{\pi^2}  \iint\limits_{[-1,1]^2\backslash J_\epsilon} \left| \frac{g_N(x)- g_N(y)}{x-y} \right|^2  \frac{dxdy}{\sqrt{1-x^2}\sqrt{1-y^2}}  \ .
 \end{align} 
 By a change of variables,
 \begin{equation*} \mathrm{I}_1(f ; N, \epsilon)
 = \frac{1}{\pi^2}  \iint_{\mathcal{B}_N} \left| \frac{f(u)- f(v)}{u-v} \right|^2 \frac{dudv}{\sqrt{1-(x_0+N^{-\alpha}u)^2}\sqrt{1-(x_0+N^{-\alpha}v)^2}}  \ ,
 \end{equation*} 
where $\mathcal{B}_N=\big[ N^\alpha(-1 + \epsilon -x_0), N^\alpha(1-\epsilon-x_0)\big]^2$. Since $f \in H^{1/2}(\R)$, by the dominated convergence theorem, we obtain for any  $0<\epsilon<1$,
 \begin{equation*} \lim_{N\to\infty} \mathrm{I}_1(f ; N, \epsilon)
 = \frac{1}{\pi^2}  \iint_{\R^2} \left| \frac{f(u)- f(v)}{u-v} \right|^2 dudv = 4 \|f\|_{H^{1/2}}^2 \ . \end{equation*} 
On the other hand, using the estimate (\ref{variance_10}),  for all sufficiently large $N$ and  for any  $0<\epsilon<\delta$,  we have
 \begin{equation*} \mathrm{I}_2(f ; N, \epsilon)
\le \frac{\mathrm{L}^2}{\pi^2\delta^2} \iint\limits_{[-1,1]^2\backslash J_\epsilon}  \frac{dxdy}{\sqrt{1-x^2}\sqrt{1-y^2}}  
=  \frac{\mathrm{L}^2}{4\delta^2} \Theta(\epsilon)^2 \ , \end{equation*}
and then, according to formula (\ref{variance_12}),  we obtain
\begin{equation} \label{variance_13}
\tilde{\Sigma}(g_N)^2= 4 \|f\|_{H^{1/2}}^2 +   \underset{N\to\infty}{O}\big(\Theta(\epsilon)\big)  \ .
\end{equation}
Finally, if we combine the estimates (\ref{variance_11}) and (\ref{variance_13}), there exists a constant $C>0$ so that  for any  $0<\epsilon<\delta$, 
\begin{equation} 
\limsup{N\to\infty}  \Var\left[ \Xi_N g_N \right]  \le 32  \|f\|_{H^{1/2}}^2 + C \Theta(\epsilon)  \ .
 \end{equation}
Since this holds for any $0<\epsilon<\delta$ and ${\displaystyle \lim_{\epsilon\searrow0} } \Theta(\epsilon) =0$, this implies formula (\ref{variance_meso}). \qed\\

\subsection{Varying exponential weights} \label{sect:UIE}


In this section, we consider unitary ensembles with varying weight of the form $\omega(x)= e^{-N V(x)}$ where the potential $V:\R\to\R$ is real-analytic.
Like in theorem~\ref{thm:UIE}, we denote this probability measure by $\mathbb{P}^V_N$.
 The condition (\ref{potential}) guarantees that (\ref{weight}) holds and the equilibrium density  $\varrho_V$ exists, see (\ref{equilibrium}).    In the following, according to (\ref{bulk}), respectively (\ref{F}), we denote
 \begin{equation*} \label{bulk_2}
  J_V = J_{\varrho_V} 
  \hspace{1cm}\text{and} \hspace{1cm}
  F_V = F_{\varrho_V}   \ .
\end{equation*}
 Moreover, by $(\ref{density_OP} - \ref{kernel_CD})$,  the spectrum $\Lambda$ of a random matrix sampled according to $\mathbb{P}^V_N$ is a determinantal process with correlation kernel
\begin{equation}\label{kernel_V}
K_N^V(x, y) = \gamma_{N-1}^2\frac{ \pi_N(x) \pi_{N-1}(y)- \pi_{N-1}(x) \pi_{N}(y)}{x-y} e^{-N\frac{V(x)+V(y)}{2}}  \ .
 \end{equation}

In the physics literature, $F_V$ is known as the {\it integrated density of states}. The set  $J_V$ corresponds to the {\it bulk} of the spectrum $\Lambda$,  it is composed of finitely many bounded open intervals and  
the equilibrium density $\varrho_V$  is smooth on  $J_V$; see \cite{Deift_al_99_b} for further references. 
%
%
One of the fundamental results of~\cite{Deift_al_99_b} is the following local asymptotics for the correlation kernel of the eigenvalue process . 
\begin{lemma}[Lem.~6.1, \cite{Deift_al_99_b}] \label{thm:sine_1}
Under the assumptions of theorem~$\ref{thm:UIE}$,
\begin{equation} \label{sine_1}
\frac{1}{N\varrho_V (x_0)} K^V_N\left(x_0 + \frac{\xi}{N\varrho_V (x_0)} , x_0 + \frac{\zeta}{N\varrho_V (x_0)} \right)= \frac{\sin \big[ \pi (\xi-\zeta)\big]}{\pi(\xi-\zeta)} + \underset{N\to\infty}{O}(N^{-1}) \ ,
\end{equation}
where the error is uniform for $x_0$ in compact subsets of $J_V$ and for $\xi,\zeta$ in compact sets of $\R$.
\end{lemma}

Actually, the non-linear steepest descent  analysis of~\cite{Deift_al_99_b}  is valid at any scales and their results implies the following sine-kernel asymptotics at mesoscopic scales. 

\begin{proposition} \label{thm:sine_V}
Under the assumptions of theorem~$\ref{thm:UIE}$, for any $\alpha\in(0,1]$,
 \begin{equation}\label{sine_0}
\frac{1}{N^\alpha} K^V_N\left(x_0 + \frac{\xi}{N^\alpha} , x_0 + \frac{\zeta}{N^\alpha} \right)= \frac{\sin \pi N \big( ( F_V (x_0 + \xi N^{-\alpha} ) - F_V(x_0 + \zeta N^{-\alpha} )\big)}{\pi(\xi-\zeta)} + \underset{N\to\infty}{O}(N^{-\alpha}) \ , 
 \end{equation}
where the error is uniform for $x_0$ in compact subsets of $J_V$ and for $\xi,\zeta$ in compact sets of $\R$.\end{proposition}

Proposition~\ref{thm:sine_V} is not formulated in~\cite{Deift_al_99_b} because the authors were interested in  universality of the local correlations and not in mesoscopic statistics. However, the proof of proposition~\ref{thm:sine_V} is a straightforward adaptation of that of lemma~\ref{thm:sine_1} and we will just review the main steps for completeness. 
First, note that $0$ is an arbitrary reference point in the definition of $F_V$. In particular, one shall  interpret the r.h.s$.$ of (\ref{sine_0}) according to (\ref{equilibrium}), namely for any $x<y$,  
$$F_V(x) - F_V(y) \simeq \frac{\# \text{ eigenvalues in }[x,y]}{N} \ .$$ 
Moreover, since the density $\varrho_V$ is smooth on $J_V$, we have 
 \begin{equation} \label{F_2}
F_V(x_0 + \xi/N) - F_V(x_0 +\zeta/N) = \varrho_V(x_0)(\xi-\zeta)N^{-1} + \underset{N\to\infty}{O}(N^{-2}) \ , 
\end{equation}
and the asymptotics (\ref{sine_1}) is a special case of (\ref{sine_0}) when $\alpha=1$.

\clearpage 

\noindent{\it Proof of proposition~\ref{thm:sine_V}.} 
We will use the Riemann-Hilbert  formulation of  \cite{Deift_al_99_b} and the formulae referenced $\{\#\}$ come from therein. We let $I=(b,a)$ be the component of $J_V$ which contains $x_0$ and for all $x\in I$,
\begin{equation}\label{F_1}
\phi(x) = \int^{a}_x \varrho_V(s)ds \ ,
\end{equation} 
see formula $\{6.7\}$ (note that in   \cite{Deift_al_99_b}, the equilibrium density is denoted by $\Psi$ instead of $\varrho_V$, cf.~$\{1.6\}$).
By $\{2.2\}$, we can write the correlation kernel 
\begin{equation} \label{CD_1}
  K_N^V(x,y) = -e^{-N(V(x)+V(y))/2 } \frac{Y_{11}(x)Y_{21}(y)- Y_{11}(y)Y_{21}(x)}{2\pi i (x-y)}\ ,  
  \end{equation}
where the $2\times2$ matrix $Y$ is the solution of a appropriate Riemann-Hilbert problem.
Transforming the  problem, cf$.$ $\{6.8-6.9\}$, the authors proved that for any $x\in I$,
\begin{equation} \label{Y_matrix}
\begin{cases} Y_{11}(x) = M_{11}(x) \exp \big[ N(V(x)+\ell +2 \pi i \phi(x))/2 \big] +  M_{12}(x) \exp \big[ N(V(x)+\ell -2 \pi i \phi(x))/2 \big] \\
Y_{21}(x) = M_{21}(x) \exp \big[ N(V(x)-\ell +2 \pi i \phi(x))/2 \big] +  M_{22}(x) \exp \big[ N(V(x)-\ell -2 \pi i \phi(x))/2 \big] 
\end{cases} ,
\end{equation}
where the $2\times2$ matrices $M(z)$ and $\frac{d}{dz} M(z)$ are uniformly bounded for all $z$ in a complex neighborhood of any point $x_0\in J_V$ and for all $N>C$, cf$.$ $\{ 5.161\}$.
Using formulae (\ref{Y_matrix}),  a little of  algebra shows that for all $x,y \in I$,
\begin{align*} 
&e^{-n (V(x)+V(y))/2 } \big(Y_{11}(x)Y_{21}(y)- Y_{11}(y)Y_{21}(x) \big) \\
&\hspace{1cm}=  e^{ i \pi n (\phi(x)-\phi(y)) } \left\{ \det M(x) -   M_{11}(x)\big(  M_{22}(x)- M_{22}(y) \big) + M_{21}(x)\big(  M_{12}(x)- M_{12}(y) \big) \right\} \\
&\hspace{1.3cm}+ e^{ -i \pi n (\phi(x)-\phi(y)) }  \left\{- \det M(x) +  M_{22}(x)\big(  M_{11}(x)- M_{11}(y) \big) - M_{12}(x)\big(  M_{21}(x)- M_{21}(y) \big) \right\}\\
&\hspace{1.3cm}+ e^{ i \pi n (\phi(x)+\phi(y)) }  \left\{  M_{21}(y)\big( M_{11}(x) -M_{11}(y) \big) - M_{11}(y)\big(  M_{21}(x) -M_{21}(y)\big) \right\}   \\
&\hspace{1.3cm}+ e^{-i \pi n (\phi(x)+\phi(y)) }  \left\{  M_{22}(y)\big( M_{12}(x) -M_{12}(y) \big) - M_{12}(y)\big(  M_{22}(x) -M_{22}(y)\big) \right\}  \\
&\hspace{1cm}= 2 i \det M(x) \sin\big(  \pi n (\phi(x)-\phi(y)) \big) + O(x-y)  \ . 
  \end{align*}
Hence, since $\det M (z) =1$ for all $z\in \C$,  by  formula (\ref{CD_1}), we obtain
\begin{equation}  \label{sine_global}
  K_N^V(x,y)= - \frac{\sin \pi N \big( (\phi (x) -\phi(y)\big)}{\pi(x-y)} + O(1) \  .
\end{equation}
Note that formula (\ref{sine_global}) holds for any points $x, y$ in the connected component $I \subseteq J_V$ which contains  the point $x_0$. To conclude it remains to observe that,  for any $L>0$ and $\alpha \in (0,1]$, if the parameter $N$ is sufficiently large, then $[x_0 -L/N^{\alpha},x_0 +L/N^{\alpha}] \subset I $ and by (\ref{F_1}), \begin{equation*}
\phi(y)-\phi(x) = \int_y^x \varrho_V(s)ds = F_V(x)-F_V(y) \ .
\end{equation*}
Hence, if we take $x= x_0 + \xi/N^\alpha$ and $y= x_0 + \zeta/N^\alpha$ with $\xi,\zeta \in [-L,L]$ in  (\ref{sine_global}) and rescale by $N^\alpha$, we obtain formula~(\ref{sine_0}).\qed \\

The correlation kernel $K_N^V$ has rank $N$ and  proposition~\ref{thm:sine_V} shows that it satisfies (\ref{sine}) for any $\alpha\in(0,1]$ with $\rho=\varrho_V$. Hence 
 theorem~\ref{thm:UIE} is a direct consequence of theorem~\ref{thm:universality}. Furthermore, in the one-cut case, the asymptotics (\ref{semiclassical_1} - \ref{semiclassical_3}) hold and we can use   
proposition~\ref{thm:variance_meso} to extend the validity of  theorem~\ref{thm:UIE} to all test functions  $f\in H^{1/2}(\R)$ which satisfies the condition (\ref{L'}). In particular, theorem~\ref{thm:clt_1} applies to the GUE and, in general, when the potential $V(x)$ is strictly convex on $\R$.  

\begin{theorem} \label{thm:clt_1}
 Let $V:\R\to\R$ be real-analytic function which satisfies the assumptions $(\ref{potential})$ and~$(\ref{1cut})$. For the eigenvalue process  of the ensemble $\mathbb{P}^V_N$, the CLT $(\ref{clt_1})$ holds for any $x_0\in J_V$, any  $0<\alpha<1$, and for all $f\in H^{1/2}\cap L^\infty(\R)$ with compact support.
 \end{theorem}

\proof  
If $X$ and $Y$ are two random variables with mean zero, by Chebychev's inequality, for any $\xi\in\R$,
\begin{equation}\label{tightness}
\left| \E{ e^{i\xi X} - e^{i\xi Y}} \right| \le  4|\xi| \sqrt{\Var\big[  X-Y\big]} \ .
\end{equation}

According to (\ref{pp}), we let  $\overline{\Xi_N^{x_0,\alpha}}f =  \Xi_N^{x_0,\alpha}f -\E{ \Xi_N^{x_0,\alpha} f}$ and $\vartheta_N(\xi; f)= \E{e^{i \xi \overline{\Xi_N^{x_0,\alpha}}f}}$ be the characteristic function of the centered linear statistics $ \Xi_N^{x_0,\alpha}f$. Recall that $\mathfrak{G}$ denotes the Gaussian field
 indexed by the Hilbert space $H^{1/2}(\R)$ and we let 
 $$ \vartheta(\xi; f)= \E{e^{i\xi \mathfrak{G}(f)}} = e^{-\frac{1}{2} \xi^2 \|f\|_{H^{1/2}}^2} \ , 
 \hspace{.6cm} \forall f\in H^{1/2}(\R) \ ,   \hspace{.5cm} \forall \xi\in\R \ .$$
 
  By the triangle inequality, for any functions $f,g \in H^{1/2}(\R)$ and for any $\xi\in\R$,
 $$ \big| \vartheta_N(\xi; f)- \vartheta(\xi; f) \big| \le
  \big| \vartheta_N(\xi; g)- \vartheta(\xi; g) \big| +  \big| \vartheta_N(\xi; f)- \vartheta_N(\xi; g) \big|
  +  \big| \vartheta(\xi; f)- \vartheta(\xi; g) \big| $$
 Furthermore, if $g\in C_0^1(\R)$, by theorem~\ref{thm:UIE}, we have
$$ \lim_{N\to\infty}   \big| \vartheta_N(\xi; g)- \vartheta(\xi; g) \big| =0 $$
 Using the estimate (\ref{tightness}) twice, since  both processes $\Xi_N^{x_0,\alpha}$ and $\mathfrak{G}$ are linear, this implies that
\begin{equation*}
\limsup{N\to\infty}  \big| \vartheta_N(\xi; f)- \vartheta(\xi; f) \big|\le
4|\xi| \left(\limsup{N\to\infty} \sqrt{\Var\big[  \Xi_N^{x_0,\alpha}(f-g)\big]}   +  \sqrt{\Var\big[ \mathfrak{G}(f-g)\big]}  \right) \ .
\end{equation*}
 If $f$ satisfies the condition (\ref{L'}), since $g\in C_0^1(\R)$, by the triangle inequality, the function $f-g$ also satisfies (\ref{L'})  and by proposition~\ref{thm:variance_meso},
$$ \limsup{N\to\infty}\Var\big[  \Xi_N^{x_0,\alpha}(f-g)\big] \le 32 \|f-g\|^2_{H^{1/2}} \ . $$ 
On the other hand, by formula (\ref{noise}), $  \sqrt{\Var[ \mathfrak{G}(f-g)]} =  \|f-g\|_{H^{1/2}}$
 and we obtain
\begin{equation} \label{variance_14} 
\limsup{N\to\infty}  \big| \vartheta_N(\xi; f)- \vartheta(\xi; f) \big|
 \le 28|\xi| \|f-g\|_{H^{1/2}} \ . 
\end{equation}
Since $C^\infty_0(\R)$ is dense in $H^{1/2}\cap L^2(\R)$, see~\cite[Thm~7.14]{Lieb_Loss},
the r.h.s.~of (\ref{variance_14})  is arbitrary small by choosing  $g\in C^1_0(\R)$ appropriately and we conclude that $\overline{\Xi_N^{x_0,\alpha}}(f) \Rightarrow \mathfrak{G}(f)$ as $N\to\infty$. \qed\\


\begin{remark} \label{rk:Venker_al}
When $V$ is strictly convex, the asymptotics (\ref{sine_0}) has also been derived in~\cite{KSSV_14} with an error which is also uniform for all potentials in a  neighborhood of $V$ (see the proof of theorem~1.7 therein).
Their method is also inspired from the results of~\cite{Deift_al_99_b} and it also applies to the slightly modified family of random matrix ensembles
\begin{equation*} \label{UIE'}
 d\mathbb{P}^{V,\mathrm{J}}_N =  \Z_{V,\mathrm{J},N}^{-1}e^{- N \tr  V(H)} \1_{\mathrm J}(H) dH \ , 
 \end{equation*}
 where $V$ is analytic and strictly convex on an interval  $\mathrm{J} \subset \R$. Hence, the results of~$\cite{KSSV_14}$ imply that theorem~\ref{thm:clt_1} holds for the  ensembles $ d\mathbb{P}^{V,\mathrm{J}}_N $ as-well.
\end{remark}

\subsection{Modified Jacobi Ensembles} \label{sect:JE}


In this section, we look at another instance of unitary ensembles given by the weight
\begin{equation} \label{weight_J}
 \omega(x) = \begin{cases} h(x) (1-x)^{\gamma_-}(1+x)^{\gamma_+} & \text{if } |x| \le1 \\
0 &\text{else}
\end{cases} \ ,
\end{equation}
where $\gamma_+, \gamma_- >-1$ and $h(x)$ is a function which is real-analytic and strictly positive on the interval $(-1-\epsilon, 1+\epsilon)$ for some $\epsilon>0$.
In this case, the probability measure  (\ref{UE}) can be written as
\begin{equation} \label{JE}
 d\mathbb{P}^\omega_N =  \Z_{\omega,N}^{-1}
  \det\big[ \omega(H) \big]\1_{\|H\| \le 1} dH \ ,
   \end{equation}
where $\|H\|$ denotes the operator norm of $H$.  We also assume  that $\omega$ is a probability density. The measure $ \mathbb{P}^\omega_N$ induces 
 a determinantal process on the eigenvalues of $H$ with correlation kernel (\ref{kernel_CD}).
 In particular, if the function $h$ is constant, the OPs with respect to $\omega$ are the classical Jacobi polynomials and their asymptotics is well-known, cf$.$~\cite{Szego} Theorem~8.21.8 and also Theorem~12.1.4. In general, the probability measure $\mathbb{P}^\omega_N$ is called the modified Jacobi unitary ensemble and the goal of this section is to derive the following asymptotics.

\begin{proposition} \label{thm:sine_J}
For any $\epsilon>0$ and $\alpha \in (0,1]$,  the correlation kernel of the modified Jacobi ensembles $\mathbb{P}^\omega_N$ with weight  $(\ref{weight_J})$ satisfies
\begin{equation} \label{sine_J}
\frac{1}{N^{\alpha}} K_N^\omega\left(x_0+ \frac{\zeta}{N^{\alpha}}, x_0+ \frac{\xi}{N^{\alpha}}\right) 
 =   \frac{\sin N \pi \big( F_\varrho(x_0+\xi N^{-\alpha}) - F_\varrho(x_0+ \zeta N^{-\alpha})\big) }{\pi ( \xi-\zeta) }
 +\underset{N\to\infty}{O} \left( N^{-\alpha}\right) \ ,
 \end{equation} 
uniformly for all $|x_0|<1-\epsilon$ and  all $\xi,\zeta$ in compact subsets of $\R$, where
\begin{equation} \label{arcsin_0}
\varrho(x)=\frac{1}{\pi \sqrt{1-x^2}} \1_{|x| \le 1} \ .
\end{equation}
\end{proposition}

The probability measure $\varrho(x)dx$ on $\R$  is called the {\bf arcsin measure} since its distribution function is given by 
 \begin{equation} \label{arcsin_1} \
F_\varrho(x)=
\begin{cases}
 \displaystyle\frac{\arcsin x}{\pi} &\text{if}\ |x|\le 1  \vspace{.2cm}\\
 \displaystyle\frac{\text{sign}\ x }{2} &\text{if}\ |x|> 1
\end{cases}\ .
\end{equation}

Proposition~\ref{thm:sine_J} implies that  $\varrho$ is the equilibrium density for the eigenvalue process  of the modified Jacobi ensembles. In contrast to the varying weights $e^{-NV(x)}$ analyzed in section~\ref{sect:UIE}, the global eigenvalues distribution is independent of the parameters of the model and it also turns out that the asymptotics of the OPs is simpler. 
Formula (\ref{sine_J}) can be deduced from the Riemann-Hilbert analysis of \cite{KV_02} by adapting the proof of proposition~\ref{thm:sine_V}. However, we will give a slightly different proof based on formula (\ref{semiclassical_1}) and the fact that the integrated density of states for the modified Jacobi ensembles is the arcsin distribution.
%
%
First, it is interesting to look at an example where we can derive proposition~\ref{thm:sine_J} using only elementary trigonometry. When $\gamma_+=\gamma_- =1/2$ and $h=1/\pi$, we denote the weight by  $\omega_0(x)= \sqrt{1-x^2}/\pi$, and the OPs which appear in the correlation kernel (\ref{kernel_CD})  are the Chebychev polynomials of the second kind. With the convention (\ref{bulk_2}), they  satisfy for all $k\ge 0$ and $x\in[-1,1]$,
\begin{equation} \label{Chebychev}
u_k(x)=  \frac{\sin\big[(k+1) \arccos x \big]}{2^k\sqrt{1-x^2}} 
\hspace{1cm}\text{and}\hspace{1cm} \gamma_k = 2^k \sqrt{2}   \ .
\end{equation}
In particular, the correlation kernel of the Chebychev process is given explicitly by 
\begin{equation}\label{Chebychev_K}
K_N^{\omega_0}(x,y) = \frac{ \sin[(N+1) \arccos x]  \sin[N \arccos y] -\sin[(N+1) \arccos y]  \sin[N \arccos x] }{\pi(1-x^2)^{1/4}(1-y^2)^{1/4} (x-y)} \ .
\end{equation}

We will need the following lemma.

\begin{lemma} \label{thm:free_K}
Let $\Psi_N$ be a function which depends on a parameter $N>0$. We define for all $|x|, |y| <1$,
\begin{equation} \label{K_psi}
 K_{\Psi_N}(x,y) = \frac{ \cos\big[ \Psi_N(x)\big] \cos\big[ \Psi_N(y)- \arccos(y) \big]
- \cos\big[ \Psi_N(y)\big] \cos\big[ \Psi_N(x)- \arccos(x) \big] }{\pi(1-x^2)^{1/4}(1-y^2)^{1/4} (x-y)} \ . 
 \end{equation}
 For any $\epsilon>0$, we have  for all  $|x|, |y| <1-\epsilon$,
  \begin{equation} \label{psi}
K_{\Psi_N}(x,y)=
 \frac{ \sin \big[ \Psi_N(y)-\Psi_N(x) \big]}{\pi(x-y)}+ \underset{|x-y|\to0}{O_\epsilon (1)} \ ,
\end{equation}
where the error term is uniform and independent of $N$.
  \end{lemma}

\proof Using the trigonometric identity
$$ \cos\big[ \Psi_N(x)- \arccos(x) \big]  = x\cos\big[ \Psi_N(x)\big] +\sqrt{1-x^2} \sin\big[ \Psi_N(x)\big] $$
we obtain, for all $|x|,|y|<1$,
\begin{align*}
& K_{\Psi_N}(x,y) \\
& = \frac{\sqrt{1-y^2}  \cos\big[ \Psi_N(x)\big] \sin\big[ \Psi_N(y)\big]
- \sqrt{1-x^2} \cos\big[ \Psi_N(y)\big] \sin\big[ \Psi_N(x) \big] }{\pi(1-x^2)^{1/4}(1-y^2)^{1/4} (x-y)} - \cos\big[ \Psi_N(x)\big] \cos\big[ \Psi_N(y)\big]  \ . 
 \end{align*}
Then, the estimate 
 $$ \left| \left( \frac{1-x^2}{1-y^2} \right)^{1/4} -1 \right| \le \frac{|x-y|}{1-y^2} $$ 
implies that for all $|x|,|y|<1-\epsilon$, 
 \begin{equation*}
 K_{\Psi_N}(x,y) = \frac{  \cos\big[ \Psi_N(x)\big] \sin\big[ \Psi_N(y)\big]
- \cos\big[ \Psi_N(y)\big] \sin\big[ \Psi_N(x) \big] }{\pi(x-y)}  + \underset{|x-y|\to0}{O_\epsilon (1)} \ ,
 \end{equation*}
 and  formula (\ref{psi}) follows from another trigonometric identity. \qed\\

The connection with the Chebychev process  is that, by (\ref{Chebychev_K}), we have $K_N^{\omega_0} =  K_{\Psi^0_N}$ with the phase  
$$\Psi^0_N(x)=(N+1) \arccos x -\pi/2 \ . $$
In particular, by (\ref{arcsin_1}), for any $x, y\in[-1,1]$, 
\begin{align}
\Psi^0_N(y) - \Psi^0_N(x) 
&\notag = N ( \arccos y  -  \arccos x)  + O(x-y) \\
&\label{arcsin_2}
= N \pi \big( F_\varrho(x) - F_\varrho(y) \big) + O(x-y)
 \end{align}
 and, according to lemma~\ref{thm:free_K}, we obtain 
\begin{equation} \label{KJ_0}
 K_N^{\omega_0}(x, y) =  \frac{ \sin \big[ N \pi \big( F_\varrho(x)-F_\varrho(y) \big) \big]}{\pi(x-y)}+ \underset{|x-y|\to0}{O (1)} \ ,
\end{equation}
uniformly for all $|x|,|y|<1-\epsilon$.


\begin{remark} \label{rk:BD} Formula $(\ref{KJ_0})$ easily yields the asymptotics of lemma~$\ref{thm:sine_J}$ and, by theorem~$\ref{thm:universality}$, this establishes the central limit theorem $(\ref{clt})$ for the Chebychev eigenvalue process . Then, universality for the modified Jacobi ensembles can be deduced from theorem~1.2 in~\cite{BD_15} since the asymptotics of the recurrence coefficients is the same for any weights of the form (\ref{weight_J}).  
Based on this approach, theorem~\ref{thm:JE} below was first proved in~\cite{BD_15} for $C^1$ test functions with compact support, cf.~theorem~1.1 therein.
In the Chebychev case, instead of using the asymptotics (\ref{KJ_0}), the authors used that the Laplace transform of the random variables $\Xi_N^{x_0,\alpha}f$ is given by a Toeplitz determinant and  computed its limit using the Strong Szeg\H{o} theorem.
\end{remark}

\begin{theorem}\label{thm:JE}
 If $(\lambda_1, \dots, \lambda_N)$ denote  the eigenvalues of a random matrix distributed according to  $\mathbb{P}^\omega_N$, $(\ref{JE})$, then for any $x_0\in (-1,1)$, any  $0<\alpha<1$, and for all $f\in H^{1/2}\cap L^\infty(\R)$ with compact support, we have as $N\to\infty$, 
\begin{equation*}
 \sum_{k=1}^N f\big( N^{\alpha}(\lambda_k-x_0)\big) - \mathbb{E}^{\omega}_N\left[  \sum_{k=1}^N f\big( N^{\alpha}(\lambda_k-x_0)\big) \right] \ \Rightarrow \ \No\big( 0, \|f\|^2_{H^{1/2}} \big) \ .
\end{equation*}
\end{theorem}


 In the remainder of this section, we will give a proof of proposition~\ref{thm:sine_J} which is inspired by the Chebychev case and lemma~\ref{thm:free_K}. The main observation is that, by theorem~\ref{thm:Kuijlaars} below, the OPs with respect to the weight (\ref{weight_J}) behave very much like the Chebychev's polynomials when $N$ is large.
By theorem~\ref{thm:universality},  proposition~\ref{thm:sine_J} implies the CLT for test functions in $C^1_0(\R)$. Moreover, since the asymptotic formulae (\ref{semiclassical_1} - \ref{semiclassical_3}) holds for the modified Jacobi ensembles, proposition~\ref{thm:variance_meso} allows us to extend the CLT for any function  $f\in H^{1/2}(\R)$ which satisfies the condition (\ref{L'}). The argument is identical to the proof of theorem~\ref{thm:clt_1}.
 The asymptotics of the OPs for the modified Jacobi ensembles has been derived using the Riemann-Hilbert method in~\cite{Kuijlaars_al_04}. In particular, we will need the following~results.

\begin{theorem}[Thm.~1.6, Thm.~1.12, \cite{Kuijlaars_al_04}] \label{thm:Kuijlaars}  
For any $\gamma_+, \gamma_- >-1$ and any function $h(x)$ which is real-analytic and strictly positive on $(-1-\epsilon, 1+\epsilon)$ for some $\epsilon>0$,  there exists $D_\infty>0$ and  $\psi  \in C^1( -1,1) $ such that the OPs with respect to $\omega(x) dx$ satisfy
$$ \pi_N(x)= \frac{D_\infty}{2^N \sqrt{\pi\omega(x)\sqrt{1-x^2}}} 
\cos\big[ (N+1/2)\arccos(x) +\psi(x)-\pi/4 \big] + \underset{N\to\infty}{O}(N^{-1}) \ , $$
uniformly for all $x$ in compact subsets of $(-1,1)$, and
$$ \frac{\gamma_N}{2^N} = \frac{\sqrt{2}}{D_\infty} \left\{ 1+ \underset{n\to\infty}{O}(N^{-1})  \right\} \ .$$ 
\end{theorem}

In a follow-up paper, \cite{KV_02}, the sine-kernel asymptotics was also derived at the local scale. For any $\epsilon, L>0$,
\begin{equation} \label{sine_J_1}
\frac{1}{N}K^\omega_N\left( x_0+\frac{\xi}{N}, x_0+\frac{\zeta}{N} \right)= \frac{\sin{[\pi \varrho( x_0)(\xi-\zeta)]}}{\pi(\xi-\zeta)}
  + \underset{N\to\infty}{O}\left(N^{-1}\right) \ ,
  \end{equation}
uniformly for all $x_0 \in [-1+\epsilon,1-\epsilon]$ and $\xi,\zeta \in [-L,L]$.  Based on the results of theorem~\ref{thm:Kuijlaars}, we obtain a first version of formula (\ref{sine_J}) which is valid as long as $|\xi-\zeta| \ge N^{-1+\epsilon}$ for any $\epsilon>0$.  Then, using local universality, we can make this asymptotics uniform for all $\xi,\zeta$ in any compact subsets of~$\R$. Hence, by combining lemmas~\ref{thm:sine_meso} and~\ref{thm:sine_micro} below, this completes the proof of proposition~\ref{thm:sine_J}.


\begin{lemma}\label{thm:sine_meso}
 For any $x_0 \in (-1,1)$, $C>0$ and $L>0$,  we have
$$ \frac{1}{N^{\alpha}} K_N^\omega\left(x_0+ \frac{\zeta}{N^{\alpha}}, x_0+ \frac{\xi}{N^{\alpha}}\right)  =   \frac{\sin N \pi \big( F_\varrho(x+\xi N^{-\alpha}) - F_\varrho(x+ \zeta N^{-\alpha})\big) }{\pi ( \xi-\zeta) }+ \underset{n\to\infty}{O}\left( N^{-\alpha} \right) \ , $$
uniformly over all  $\xi,\zeta \in[-L,L]$ such that  $|\xi-\zeta| \ge C N^{-1+\alpha}$.
\end{lemma}

\proof For any $x\in(-1,1)$ and $N\ge 1$, we let
\begin{equation*} \Psi_N(x) = (N +1/2)\arccos(x) +\psi(x)-\pi/4 \ .
 \end{equation*}
Formula (\ref{kernel_CD}) and the asymptotics of theorem~\ref{thm:Kuijlaars} implies that for any $x,y \in [-1+\epsilon, 1-\epsilon]$, 
\begin{equation*}
 K_N^\omega(x, y) = K_{\Psi_N}(x,y) + O\left(  \frac{1}{|x-y|N} \right) \ ,
 \end{equation*}
where $ K_{\Psi_N}$ is given by formula (\ref{K_psi}). Hence, by lemma~\ref{thm:free_K},
\begin{equation*}
 K_N^\omega(x, y) =\frac{ \sin \big[ \Psi_N(y)-\Psi_N(x) \big]}{\pi(x-y)}+O\left(1+  \frac{1}{|x-y|N} \right) \ .
 \end{equation*}
Since $\psi\in C^1(-1,1) $, like in the Chebychev case, we have
$$ \Psi_N(y) - \Psi_N(x)  = N \pi \big( F_\varrho(x) - F_\varrho(y) \big) + O(x-y) \ , $$
and 
\begin{equation} \label{KJ_1}
 K_N^\omega(x, y) =\frac{ \sin \big[N \pi \big( F_\varrho(x) - F_\varrho(y) \big)\big]}{\pi(x-y)}+O\left(1+  \frac{1}{|x-y|N} \right) \ .
 \end{equation}
To conclude it remains to take $x=x_0 + \xi  N^{-\alpha}$, $y=x_0 + \zeta  N^{-\alpha}$ and rescale by $N^{-\alpha}$.  In this regime, the error in formula (\ref{KJ_1}) is of order $N^{-\alpha} + 1/|\xi-\zeta|N$ and lemma~\ref{thm:sine_meso} follows immediately. \qed\\


A fundamental observation due to K. Johansson is that, in the regime $|\xi-\zeta| \le N^{-1+\alpha}$, the difference  $N^{-\alpha}(\xi-\zeta)$ is {\it microscopic} and we can deduce the asymptotics of the kernel using only local universality  considerations.

\begin{lemma} \label{thm:sine_micro}
 For any $x_0 \in (-1,1)$, $C>0$ and $L>0$,   
$$ \frac{1}{N^{\alpha}} K_N^\omega\left(x_0+ \frac{\zeta}{N^{\alpha}}, x_0+ \frac{\xi}{N^{\alpha}}\right) =   \frac{\sin N \pi \big( F_\varrho(x+\xi N^{-\alpha}) - F_\varrho(x+ \zeta N^{-\alpha})\big) }{\pi ( \xi-\zeta) }+ \underset{n\to\infty}{O}\left( N^{-\alpha} \right) \ , $$
uniformly over all  $\xi,\zeta \in[-L,L]$ such that  $|\xi-\zeta| \le C N^{1-\alpha}$.
\end{lemma}

\proof  We define
\begin{equation}\label{micro}
 \tilde x_0 = x_0 + \frac{\xi+\zeta}{2} N^{-\alpha} \ ,
\hspace{.5cm} 
\tilde \xi =  \frac{\xi-\zeta}{2} N^{1-\alpha}
\hspace{.5cm} \text{and} \hspace{.5cm}
 \tilde \zeta = \frac{\zeta-\xi}{2} N^{1-\alpha} \ .
\end{equation}
 Since $\alpha>0$, when $N$ is sufficiently large, there exists $\epsilon>0$ such that $\tilde x_0 \in [-1+\epsilon,1-\epsilon]$. Moreover, the assumption $|\xi-\zeta| \le C N^{-1 +\alpha}$ implies that $\tilde \xi , \tilde \zeta \in [-C/2,C/2]$ and we can use the asymptotics (\ref{sine_J_1}), 
\begin{equation*}
\frac{1}{N} K_N^\omega\left(\tilde x_0+\frac{\tilde\xi}{N},\tilde x_0+\frac{\tilde\zeta}{N} \right)= \frac{\sin{[\pi \varrho(\tilde x_0)(\tilde\xi-\tilde\zeta)]}}{\pi(\tilde\xi-\tilde\zeta)}
  +O\left(N^{-1}\right) \ .
  \end{equation*}
By (\ref{micro}), this implies that
\begin{align}  
& \notag
 \frac{1}{N}K_N^\omega\left(x_0+ \frac{\zeta}{N^{\alpha}}, x_0+ \frac{\xi}{N^{\alpha}}\right)  
  = \frac{\sin{[ N^{1-\alpha} \pi \varrho(\tilde x_0)(\xi-\zeta)]}}{N^{1-\alpha} \pi (\xi-\zeta)} + O\left(N^{-1}\right) \ , \\
&  \label{KJ_3}
\frac{1}{N^{\alpha}} K_N^\omega\left(x_0+ \frac{\zeta}{N^{\alpha}}, x_0+ \frac{\xi}{N^{\alpha}}\right) 
=  \frac{\sin{[ N^{1-\alpha} \pi \varrho(\tilde x_0)(\xi-\zeta)]}}{\pi (\xi-\zeta)}
   + O\left(N^{-\alpha} \right)   \ .
\end{align}
On the other hand, by definition, two taylor expansions yields when $|\xi-\zeta| \le C N^{-1 +\alpha}$,
 \begin{align*} F_\varrho\left(x_0 +N^{-\alpha} \xi \right)- F_\varrho\left(x_0+N^{-\alpha} \zeta\right)  
&= \varrho\left( x_0 + \frac{\xi+\zeta}{2} N^{-\alpha} \right) N^{-\alpha} (\xi - \zeta) 
+ O\left(|N^{-\alpha} (\xi - \zeta) |^3\right)  \\
&= \varrho\left( \tilde x_0\right)  N^{-\alpha} (\xi - \zeta)  + O\left(|\xi-\zeta| N^{-2-\alpha}\right)  \ ,
\end{align*}
and we obtain
\begin{equation} \label{KJ_4}
\frac{\sin N \pi \left( F_\varrho\left(x_0 +N^{-\alpha} \xi \right)- F_\varrho\left(x_0+N^{-\alpha} \zeta\right) \right)}{\pi(\xi-\zeta)}  =  \frac{\sin{[ N^{1-\alpha} \pi \varrho(\tilde x_0)(\xi-\zeta)]}}{ \pi (\xi-\zeta)} + O\left(N^{-1-\alpha}\right)
\end{equation}
We conclude the proof by combining formulae (\ref{KJ_3}) and  (\ref{KJ_4}).
\qed \\


\begin{remark} \label{rk:strong_asymptotics}   To  prove lemma~$\ref{thm:sine_J}$  for all $\alpha \in (0,1]$, it is important to use the uniform asymptotics of theorem~$\ref{thm:Kuijlaars}$ and $(\ref{sine_J_1})$ with the optimal error of order $1/N$. 
 There are other methods than the Riemann-Hilbert steepest descent to compute the asymptotics of OPs for a non-varying measure  and  prove local universality, e.g$.$  the methods developped by Levin and Lubinsky, \cite{Lubinsky_09a,Lubinsky_09b}. However, these methods usually provide weaker asymptotics which yields the  sine-kernel only at small scales; see also remark~\ref{rk:Levin_Lubinsky} below.
\end{remark}

\section{Regularized characteristic polynomial and log-correlated Gaussian processes} \label{sect:fbm}

The goal of this section is to elaborate on the connection between logarithmically correlated Gaussian processes (1/f-noise) and random matrix theory. It was established in~\cite{HKO_01} and \cite{FKS_13} that the logarithm of the modulus of the characteristic polynomial of a CUE, respectively GUE, random matrices converge weakly to Gaussian  generalized functions (random tempered distributions) whose correlation kernels have a logarithmic singularity at 0. 
Based on the so-called {\it freezing transition scenario}, this motivates some recent conjectures for the distributions of the extreme values of these polynomials, as well as for the extreme value of  the Riemann Zeta function on some interval of the critical line, see \cite{FHK_12, FS_15} and references therein. This also suggests that the characteristic polynomials of random matrices give raise to regularizations of the so-called Gaussian Multiplicative Chaos measures introduced by Kahane, which play an important role in some recent physical theory, such as conformal field theory, disordered systems, Liouville quantum gravity, etc, \cite{DRSV_14, Webb_15}. 
%
 In the following, we consider a random Hermitian matrix $H$ distributed according  to the  unitary invariant measure $\mathbb{P}^\omega_N$, (\ref{UE}). We  will not look directly at the characteristic polynomial of the matrix $H$ but the following regularization at mesoscopic scales.
Let  $0<\alpha<1$, $x_0\in \R$, $\eta>0$, $z_t = t + i \eta$, and define
\begin{equation} \label{W1}
W_N(t) = \log\left| \det\left[ H - x_0- z_t N^{-\alpha} \right] \right| -  \log\left| \det\left[ H - x_0- z_0 N^{-\alpha} \right] \right| \ . 
\end{equation}
This object was introduced in \cite{FKS_13} and it was proved that  if $H$ is a GUE matrix, then the random process 
$t\mapsto W_N(t) - \E{W_N(t)}$ converges weakly in $L^2[a,b]$ ($a,b \in \R$) to a logarithmically correlated Gaussian process $B_0$ defined below. 

\begin{definition} \label{FBM}
 The $\eta$-regularized fractional Brownian motion with Hurst exponent $H=0$ is a real-valued Gaussian process $B_0$ characterized by the following properties:
\begin{itemize}
\item[i)]  $B_0$ is a continuous process with mean $0$ and $B_0(0)=0$ almost surely. 
\item[ii)]  $B_0$ has stationary increments.
\item[iii)]  $\displaystyle \Var\big[ B_0(t)  \big] = \frac{1}{2}\log\left(1+ \frac{t^2}{4\eta^2}  \right)$ for any $t\in \R$.
\end{itemize}  
\end{definition}  

\noindent We refer to \cite{FKS_13}  for some background and references on fractional Brownian motion. Let us just point out that the process $B_0$ has the following representation, for any $t\in\R$,
$$ B_0(t) = \Re \bigg\{ \int_0^\infty e^{-\eta s} ( e^{-i ts}-1 ) \frac{dZ_s}{\sqrt{2s}} \bigg\} \ , $$    
where $Z$ is a complex Brownian motion with unit variance.  Inspired by Riemann-Hilbert asymptotics obtained by Krasovsky in~\cite{Krasovsky_07},  the authors of~\cite{FKS_13} computed the limits of the Laplace transform of the random variable $W_N(t)$ for any $t\in\R$ and show that the finite-dimentionnal distributions of  $ W_N - \E{W_N}$  converges to that of $B_0$. In the following, we generalize this result to other unitary invariant ensembles using the central limit theorem \ref{thm:clt_1}. We suppose that the weight $\omega$ satisfies (\ref{weight}) and the one-cut condition, $J_{\varrho_\omega} = (-1,1)$, so that the estimates of section~\ref{sect:variance} hold.
 Although this condition should not be relevant, we have not derived the necessary  variance estimates in the multi-cut regime.
 First, observe that the random variable $W_N(t)$ is a linear statistic,
\begin{align}  W_N(t) 
&\notag = \Re\bigg\{ \log \det \bigg[ \frac{M - x_0- z_t N^{-\alpha}}{M - x_0- z_0 N^{-\alpha}}\bigg]\bigg\} \\
&\notag =  \Re\bigg\{ \tr\bigg[ \log \bigg(  \frac{M - x_0- z_t N^{-\alpha}}{M - x_0- z_0N^{-\alpha}}\bigg) \bigg]\bigg\} \\
&\label{W0} =  \Xi_N^{x_0,\alpha} g_t  
\end{align}  
 where the function $\displaystyle g_t(x) = \Re\bigg\{ \log \bigg(  \frac{x- z_t }{x- z_0}\bigg)\bigg\}$ is defined using the principal branch of the logarithm and $z_t = t + i \eta$. 
It is easily seen that,  even though $g_t \notin L^1(\R)$, its Fourier transform is well defined in $L^2(\R)$ and, by lemma~\ref{thm:Frullani} below, it is given by 
\begin{equation} \label{W2}
 \widehat{g_t}(u) = (1-e^{- 2\pi i u t})  \frac{e^{-2\pi |u| \eta}  }{2|u|} \ . 
 \end{equation}

\begin{lemma} \label{thm:Frullani}
 For any $\eta>0$ and $x, t\in \R$,  we have
\begin{equation*} \int_\R e^{2\pi iu x}  (1-e^{-2\pi i u t})  \frac{e^{- 2\pi\eta |u|}  }{2|u|} du 
= g_t(x) =  \frac{1}{2}\log\left( \frac{(x-t)^2+ \eta^2}{x^2+ \eta^2}  \right)  \ .
\end{equation*} 
\end{lemma}

\proof  This identity is classical and it can be proved by observing that, for any $t>0$,
$$ \frac{1- e^{- 2\pi it u}}{2|u|} = i \sgn(u) \pi  \int_0^t e^{-2\pi i s u } ds \ , $$
and, by Fubini's theorem,
\begin{align*} \int_\R e^{2\pi iu x}  (1-e^{-2\pi i u t})  \frac{e^{- 2\pi\eta |u|}  }{2|u|} du &= i  \pi  \int_{0}^t  \int_\R e^{2\pi iu (x-s)} \sgn(u) e^{- 2\pi \eta |u|} du \ ds   \\
&= -2\pi \int_{0}^t \Im\bigg\{  \int_0^\infty e^{- 2\pi u (\eta-i  (x-s))}  du \bigg\} ds \\
&= - \Im\bigg\{ \int_{0}^t  \frac{ds}{\eta - i (x-s) } \bigg\} \ . 
\end{align*} 
We conclude by observing that, by definition, for any $t>0$, 
\begin{equation} \label{W3}
 g_t(x) = \Re\bigg\{ \int_{x-t}^x \frac{dv}{v+i\eta}\bigg\} \ . 
 \end{equation}
The proof in the case  $t<0$ is almost identical. \qed\\ 

 At the end of this section, we check that the test functions $g_t$ satisfy the assumptions of proposition~\ref{thm:variance_meso} so that for any $|x_0|<1$ and $0<\alpha<1$,
 \begin{equation} \label{W4}
 \limsup{N\to\infty} \Var\left[ \Xi_N^{x_0,\alpha} g_t \right]  \le 32\ \| g_t\|_{H^{1/2}}^2 
 \end{equation}
 and we can apply the CLT (\ref{clt_1}), cf.~the proof of theorem~\ref{thm:clt_1}. Namely,  for any $t_1<\dots<t_k$ and $\xi_1,\dots, \xi_k \in \R$, letting $f = \sum_{j=1}^k  \xi_j g_{t_j}$, we obtain
\begin{equation} \label{W5}
\Xi_N^{x_0,\alpha} f - \mathbb{E}^\omega_N\left[ \Xi_N^{x_0,\alpha} f\right] \ \Rightarrow \ \No\big( 0, \|f\|^2_{H^{1/2}} \big) \ ,
\end{equation}
where
$$  \|f\|^2_{H^{1/2}}  = \sum_{l, j = 1}^k \xi_l \xi_j \langle g_{t_i} , g_{t_j} \rangle_{H^{1/2}} \ .  $$ 
 Moreover, by formula (\ref{W2}), 
$$ \widehat{g_t}(u) \widehat{g_s}(-u) |u| = (1 -e^{-i2\pi u t} -e^{i2\pi u s} + e^{i2\pi u(s-t)}  )  \frac{e^{-4\pi |u| \eta}  }{4|u|} \ , $$ 
and, according to lemma~\ref{thm:Frullani} with $x=0$, we obtain for any $t, s \in \R$
\begin{align}
\langle g_t , g_s \rangle_{H^{1/2}} 
&\notag = \int_\R \widehat{g_t}(u) \widehat{g_s}(-u) |u| \\
&\label{W6} =     \frac{1}{4} \left\{ \log\left(1+ \frac{t^2}{4\eta^2}  \right)+ \log\left(1+ \frac{s^2}{4\eta^2}  \right) -  \log\left(1+ \frac{(t-s)^2}{4\eta^2}  \right) \right\} \ .
\end{align}

\noindent Since we have established that $  \Xi_N^{x_0,\alpha} g_t = W_N(t)$,  according to definition~\ref{FBM},  formulae (\ref{W5} - \ref{W6}) imply that for any $k\in\N$,
\begin{equation} \label{W7}
 \big( W_N(t_1)-  \mathbb{E}^\omega_N\left[ W_N(t_1)\right] ,\dots, W_N(t_k)-  \mathbb{E}^\omega_N\left[ W_N(t_k)
 \right] \big)  \ \Rightarrow \   \big( B_0(t_1),\dots, B_0(t_k) \big) \ . 
 \end{equation}
Note that the fact that the Gaussian process $B_0$ has independent increments follows immediately from the covariance structure (\ref{W6}) and the continuity of its sample paths follows from Kolmogorov's theorem. Following~\cite[Thm.~2.3]{FKS_13}, the convergence (\ref{W7}) of the finite-dimensional distributions and the estimate (\ref{W4}) allows to conclude that the random process $W_N$ converges in distribution to $B_0$ in an appropriate function space. 

\begin{theorem} \label{thm:FBM} Let $\omega$ be any positive function satisfying  $(\ref{weight})$ and such that the support of its equilibrium measure satisfies $J_{\varrho_\omega} = (-1,1)$. For any $|x_0|<1$, any $0<\alpha <1$, and any $a,b \in \R$ such that $a<b$,  the stochastic process
 $ W_N-  \mathbb{E}^\omega_N\left[ W_N\right]$, $(\ref{W1})$, converges weakly as $N\to\infty$ in $L^2[a,b]$ to the  $\eta$-regularized fractional Brownian motion $B_0$ with Hurst exponent $H=0$.   
\end{theorem}

To complete the proof of theorem~\ref{thm:FBM}, it remains to check that  the functions $g_t$ satisfies the assumptions of proposition~\ref{thm:variance_meso} for all $t \in \R$.

\begin{lemma}For any $t\in \R$, the function $g_t \in H^{1/2}(\R)$ and it satisfies the hypothesis $(\ref{L'})$. 
\end{lemma}

\proof Without loss of generality we suppose that $t>0$. Formula (\ref{W3})
 implies that 
$$ g_t(x) - g_t(y) = \Re\bigg\{ \int_{ \mathscr{C}_{x,y}} \frac{dz}{z}\bigg\} \  , $$
where for any $(x,y)\in \R^2$,
$$ \mathscr{C}_{x,y} =   \begin{cases}
 \{ v + i \eta  :  v \in [y-t, y] \uplus [x-t , x] \}  &\text{if } y < x- t  \\
 \{ v + i \eta  :  v \in [y-t, x-t] \uplus [y, x] \} &\text{if } x-t<y< x \\
  \{ v + i \eta  :  v \in [x-t, y-t] \uplus [x, y] \} &\text{if } x< y <  x+ t \\ 
 \{ v + i \eta  :  v \in [x-t, x] \uplus [y-t , y] \}  &\text{if }  x<y-t  
\end{cases} \ .$$

Note that, in all four cases, the length of the contour $ |  \mathscr{C}_{x,y} | = 2 \min\{ t,  |x-y| \}$, and there exists $C_t>0$ and a 
continuous function $h:\R\to\R_+$ such that 
$$ \max_{z\in \mathscr{C}_{x,y}} |z|^{-1}  \le h(x)/ 2 
\hspace{.8cm}\text{and}\hspace{.8cm} |h(x)| \le 4/ |x| \hspace{.7cm} \text{for all } |x|\ge C_t \ .$$
This implies that
\begin{align*} \big| g_t(x) - g_t(y) \big| & \le |  \mathscr{C}_{x,y} | \max_{z\in \mathscr{C}_{x,y}} |z|^{-1}  \\
& \le \min\{ t,  |x-y| \} h(x) \ ,
\end{align*}

By Fubini's theorem, we conclude that 
  \begin{align*} 
\iint\left| \frac{g_t(x)- g_t(y)}{x-y} \right|^2 dxdy  = \int h (x)^2 dx  \int  \left(\frac{\min\{t, |x-y|\}}{x-y}\right)^2 dy  < \infty \ .
\end{align*}

Moreover,  by construction, for all $|x|\ge C_t$, we have
$$  \sup\left\{     \left| \frac{g_t(x)- g_t(y)}{x-y} \right| : |y| \le |x| \right\}  
\le \frac{4}{|x|}
$$
so that  the hypothesis (\ref{L'}) holds. \qed\\

\section{The Gaussian Unitary Ensemble} \label{sect:GUE}

The GUE  (\ref{GUE}) was introduced by E. Wigner as a model to describe scattering resonances of Heavy  nuclei and it is certainly the Hermitian matrix model which received most attention.
In particular, in addition to be unitary invariant, the entries of  a GUE matrix are independent Gaussian random variables. The GUE falls in the general class discussed in section~\ref{sect:UIE} with weight $\omega(x)= e^{-N x^2}$. Hence, theorem~\ref{thm:UIE} implies that its eigenvalue process  converges at mesoscopic scales to the $H^{1/2}$- Gaussian field $\mathfrak{G}$. 
In fact, another proof valid for Gaussian $\beta$-ensembles,   appeared  previously in~\cite{BEYY_14}.
 The goal of this section is to derive the GUE kernel asymptotics from the classical integral formulae for the Hermite polynomials rather than by solving a Riemann-Hilbert problem. 
 We proceed like in section~\ref{sect:JE}. First, in section~\ref{sect:GUE_2}, we produce the global asymptotics of the correlation kernel.  Then, in section~\ref{sect:GUE_3}, we make this asymptotics uniform by looking at the microscopic regime.

\subsection{Plancherel-Rotach asymptotics}

The first observation is that the GUE weight satisfies $\omega(x)= \omega_G(\sqrt{2N} x)$ where  $\omega_G(x)= e^{-x^2}$ does not depend on the dimension $N$. Moreover, the OPs with respect to the Gaussian weight $\omega_G$ are the classical Hermite polynomials, for all $k \ge 0$,
\begin{equation}\label{Hermite_0}
\pi_k(x)= e^{x^2}\left(-\frac{1}{2}  \frac{d}{dx}\right)^k e^{-x^2} 
\hspace{1cm}\text{and}\hspace{1cm}
\gamma_k = \sqrt{\frac{2^k}{\sqrt{\pi}k!}}  \ .
\end{equation}
If we let 
\begin{equation}\label{Hermite_0'}
\phi_k(x)=\sqrt{\omega_G(x)} \gamma_k \pi_k(x)  \ , 
\end{equation}
then, according to formula (\ref{kernel_CD}), the correlation kernel of the GUE eigenvalue process  is given by
\begin{equation} \label{K_GUE}
K_N^{\G}(x,y)= \sqrt{N} K_N^{\omega_G}(\sqrt{N} x,\sqrt{N}y) 
\end{equation}  
where
\begin{equation}  \label{K_G} 
K_N^{\omega_G}(x, y) 
= \sqrt{\frac{N}{2}} \frac{ \phi_N(x) \phi_{N-1}(y)- \phi_{N-1}(x) \phi_{N}(y)}{x-y} \ .
\end{equation}
The functions $\phi_k$ are usually called the Hermite (wave) functions and they form an orthonormal basis of $L^2(\R)$. Moreover, they have the following asymptotics.

\begin{proposition} \label{thm:Hermite}
Let, for all $|x| <1$ and $N>0$, 
\begin{equation}  \label{H}
H(x)= \arccos x - x \sqrt{1-x^2} 
\hspace{1cm}\text{and}\hspace{1cm}
\Psi_N(x) = NH(x) + \frac{\arccos x}{2} - \frac{\pi}{4} \ .
\end{equation}
There exists  two sequences of functions $\Lambda_N$ and $\tilde\Lambda_N$ which are smooth on $(-1,1)$ such that for any $\epsilon>0$ and for all $|x|\le1-\epsilon$,
\begin{align}
\phi_{N}(\sqrt{2N}x)
&\label{Hermite_1}=\frac{\eta_{N+1}}{\big(N(1-x^2) \big)^{1/4}}\bigg\{\cos \big[\Psi_N(x)\big] + \tilde\Lambda_N(x)  
+\underset{N\to\infty}{O_\epsilon}\big(N^{-2}\big) \bigg\}  \\
\phi_{N-1}(\sqrt{2N}x)
&\label{Hermite_1'}=\frac{\eta_{N}}{\big(N(1-x^2) \big)^{1/4}}\bigg\{\cos \big[\Psi_N(x)-\arccos(x)\big] + \Lambda_N(x)  
+\underset{N\to\infty}{O_\epsilon}\big(N^{-2}\big) \bigg\} \ ,
\end{align}
where
\begin{equation} \label{nu}
 \eta_N=  \sqrt{\frac{e^NN!}{\pi^{3/2}N^N}} = \frac{2^{1/4}}{\sqrt{\pi}} + \underset{N\to\infty}{O}\big(N^{-1}\big) \ .
 \end{equation}
 Moreover, there exists a universal constant $C>0$ such that for all $|x|<1$ and $N> 0$,
\begin{equation}\label{Lambda} 
|\Lambda_N(x)\vee \tilde\Lambda_N(x)| \le \frac{C}{N(1-x^2)^{3/2}}
\hspace{1cm}\text{and}\hspace{1cm}
| \Lambda'_N(x)\vee \tilde\Lambda_N'(x) | \le  \frac{C}{(1-x^2)^{5/2}} \ .  
 \end{equation}
\end{proposition}

\proof Thanks to Rodrigues' formula, (\ref{Hermite_0}), the Hermite functions have the following integral representation
\begin{equation*}
\phi_k(x)=  e^{-x^2/2} \sqrt{\frac{k!}{\sqrt{\pi}2^k}} \frac{1}{2\pi i} \oint_{|z|=1} z^{-k-1} e^{-(z-x)^2} dz  \ .
\end{equation*}
The saddle point analysis for this integral was performed in a  seminal paper by Plancherel and Rotach, \cite[ formula~7]{PR_29}. If $\eta_N$ is given by (\ref{nu}),
\begin{equation}  \label{H_0}
H(x)= \arccos x - x \sqrt{1-x^2} 
\hspace{1cm}\text{and}\hspace{1cm}
\varphi(x) = \frac{\arccos x}{2} + \frac{\pi}{4} \ ,
\end{equation}
they obtained for any $k\in \N$ and $|x|<1$,
\begin{align}\label{PR} 
\phi_{N-1}(\sqrt{2N}x)
&=\frac{ \eta_N}{N^{1/4}(1-x^2)^{1/4}}\\
&\notag\hspace{.3cm}
\times\bigg\{\sum_{ s =0}^{k-1} N^{-s} \sum_{i=0}^{2s}  C_{s,i}  \frac{\cos \big[N H(x)- (2s+2i+1) \varphi(x) \big]}{(1-x^2)^{(s+i)/2}}
 +\underset{N\to\infty}{O}\left(\frac{1}{N^k(1-x^2)^{3k/2}}\right)\bigg\} \ .
\end{align} 

It is remarkable that they managed to obtain the full asymptotic expansion of the Hermite functions. In fact, to obtain the mesoscopic asymptotics of the GUE kernel, it suffices to take $k=2$ in (\ref{PR}), then the coefficients in the expansion are  $C_{0,0}=1$,  $C_{1,0}=0$,  $C_{1,1}=3/16$ and $C_{1,2}= 5/48$ according to \cite{PR_29}.
In this case, we deduce from formula (\ref{PR}) that uniformly for all $x$ in compact subsets of $(-1,1)$,
\begin{equation}\label{Hermite_3} 
\phi_{N-1}(\sqrt{2N}x)
=\frac{ \eta_N}{N^{1/4}(1-x^2)^{1/4}} \bigg\{\cos \big[N H(x)- \varphi(x) \big] + \Lambda_N(x)
 +\underset{N\to\infty}{O}\left(N^{-2}\right)\bigg\} \ ,
\end{equation} 
and the function $\Lambda_N$ is smooth on $(-1,1)$ and  satisfies 
\begin{equation}\label{Lambda_1} 
|\Lambda_N(x)| \le \frac{1}{N(1-x^2)^{3/2}}
\hspace{1cm}\text{and}\hspace{1cm}
| \Lambda'_N(x) | \le \frac{|H'(x)|}{(1-x^2)^{3/2}} + \frac{3}{N(1-x^2)^{5/2}} \ .  
\end{equation}
According to formula (\ref{Hermite_3}), if we let $x_N =  \sqrt{\frac{N}{N+1}}  x $,  we obtain 
\begin{align}\label{Hermite_4} 
\phi_{N}(\sqrt{2N}x)
=\frac{\eta_{N+1}}{\big(N(1-x^2) +1\big)^{1/4}}\bigg\{\cos \big[(N+1) H(x_N)- \varphi(x_N) \big] + \Lambda_{N+1}(x_N)
 +\underset{N\to\infty}{O}\left(N^{-2}\right)\bigg\} \ .
\end{align} 
Moreover, by (\ref{H_0}), we see that  for any $|x|<1-\epsilon$,
\begin{align}
 (N+1) H\left( \sqrt{\frac{N}{N+1}}  x\right) - \varphi\left( \sqrt{\frac{N}{N+1}}  x\right) 
& \notag = N H(x) -\varphi(x) - \frac{xH'(x)}{2}  + H(x)  +O(N^{-2}) \\
& \label{phase}
 = N H(x) -\varphi(x) + \arccos x  +O(N^{-2}) \ .
\end{align}
This identity is remarkable because if $\Psi_N$ is defined according to (\ref{H}) 
and we substitute  (\ref{phase}) in formula (\ref{Hermite_4}), we obtain
\begin{equation*}
\phi_{N}(\sqrt{2N}x)
=\frac{\eta_{N+1}}{\big(N(1-x^2) +1\big)^{1/4}}\bigg\{\cos \big[\Psi_N(x)\big] + \Lambda_{N+1}(x_N)  +O\big(N^{-2}\big) \bigg\}  \ .
\end{equation*}
Moreover, for any $|x| \le 1-\epsilon$,
\begin{equation*}
\left(\frac{1-x^2}{1-x^2 +1/N}\right)^{1/4} = 1+ \frac{1}{N(1-x^2)} +O\big(N^{-2}\big) \ ,
\end{equation*}
and this implies that
\begin{equation}\label{Hermite_5}
\phi_{N}(\sqrt{2N}x)
=\frac{\eta_{N+1}}{N^{1/4}(1-x^2)^{1/4}}\bigg\{\cos \big[\Psi_N(x)\big] + \tilde\Lambda_N(x)  +O\big(N^{-2}\big) \bigg\}  \ ,
\end{equation}
where,
\begin{equation*}
\tilde\Lambda_N(x) = \Lambda_{N+1}(x_N) +  \frac{\cos\big[\Psi_N(x)\big]}{N(1-x^2)} \ .
\end{equation*}
  Since $x_N =  \sqrt{\frac{N}{N+1}}  x $, using the estimates  (\ref{Lambda_1}),  we see that the function $\tilde\Lambda_N$ is smooth on $(-1,1)$ and it satisfies
\begin{equation}\label{Lambda_2} 
|\tilde\Lambda_N(x)| \le \frac{2}{N(1-x^2)^{3/2}}
\hspace{1cm}\text{and}\hspace{1cm}
| \Lambda'_N(x) | \le \frac{3|H'(x)|}{(1-x^2)} + \frac{6}{N(1-x^2)^{5/2}} \ .  
\end{equation}
Finally, by (\ref{H}), the function $H\in C^1(-1,1)$ and 
\begin{equation} \label{H_1}
H'(x)= -2 \sqrt{1-x^2} \ ,
\end{equation}
so that the estimates (\ref{Lambda}) follow from (\ref{Lambda_1})  and (\ref{Lambda_2}).\qed\\

\subsection{The global asymptotics} \label{sect:GUE_2}

Proposition~\ref{thm:Hermite} encompasses most of the technical work to prove formula (\ref{sine}) for the GUE kernel. In this section, we will derive the global asymptotics of the GUE kernel
 based on the method developed in section~\ref{sect:JE}.

\begin{lemma} \label{thm:sine_meso'} 
Let, for all $|x| <\sqrt{2}$, 
 \begin{equation}\label{sc}
 \varrho_{\s}(x)= \frac{\sqrt{2-x^2}}{\pi} \1_{[-\sqrt{2},\sqrt{2}]}(x) 
 \hspace{1cm}\text{and}\hspace{1cm}
 F_{\s}(x)= \int_0^x  \varrho_{\s}(u) du  \ .
\end{equation}
For any $\epsilon>0$ and for all $|x| , |y| \le \sqrt{2} -\epsilon$, we have 
 \begin{equation} \label{gue_0}
K_N^{\G}(x,y) 
= \frac{ \sin \big[ \pi N \big( F_{\s}(x)-F_{\s}(y) \big)\big]}{\pi(x-y)}+ \underset{N\to\infty}{O_\epsilon}\left(1 + \frac{1}{|x-y|N^2}\right) \ .
\end{equation}
\end{lemma}

\proof We define
\begin{equation} \label{Psi_1}
  \tilde \Psi_N(x) = \Psi_N(x) - \arccos x \ .
\end{equation}
By formulae (\ref{H}) and (\ref{H_1}), for all  $|x|\le1-\epsilon/\sqrt{2}$,
\begin{equation} \label{Psi_2}
\big|\Psi_N(x) - \Psi_N(y)\big| \le (2N + \epsilon^{-1/2}) |x-y| \ .
\end{equation}
and the same bound holds for  $\tilde\Psi_N$. 
 By (\ref{K_G}) and proposition~\ref{thm:Hermite}, we get
\begin{align}
&\notag  \sqrt{2N}K_N^{\omega_G}(\sqrt{2N}x, \sqrt{2N} y) 
 =   \frac{ \pi \eta_N \eta_{N+1}}{\sqrt{2}} K_{\Psi_N}(x,y) \\
&\notag\hspace{.3cm}  + \frac{\eta_N \eta_{N+1}}{(1-x^2)^{1/4}(1-y^2)^{1/4}} \bigg\{ 
\frac{\tilde\Lambda_N(x)\Lambda_N(y)-\tilde\Lambda_N(y)\Lambda_N(x)}{x-y}
+\frac{\cos \big[\Psi_N(x)\big] \Lambda_N(y) - \cos \big[\Psi_N(y)\big] \Lambda_N(x)}{x-y} \\
&\label{gue_1}
\hspace{.6cm} 
+\frac{\cos \big[\tilde\Psi_N(x)\big] \tilde\Lambda_N(y)-\cos\big[\tilde\Psi_N(y)\big] \tilde\Lambda_N(x)}{x-y}
+\underset{N\to\infty}{O_\epsilon}\left(\frac{1}{|x-y|N^2}\right) \bigg\} \ , 
\end{align}
where, in the first term, the kernel is given by formula~(\ref{K_psi}). By lemma~\ref{thm:free_K},  this term yields the sine-kernel asymptotics. Indeed, by (\ref{nu}), we have  for all  $|x|\le1-\epsilon/\sqrt{2}$,
\begin{align} 
 \frac{ \pi \eta_N \eta_{N+1}}{\sqrt{2}} K_{\Psi_N}(x,y) 
 & \notag = \left\{ 1 + O\big(N^{-1}\big)  \right\}
 \frac{ \sin \big[ \Psi_N(y)-\Psi_N(x) \big]}{\pi(x-y)}+ O_\epsilon(1)  \\
 & \label{gue_4}=  \frac{ \sin \big[ \Psi_N(y)-\Psi_N(x) \big]}{\pi(x-y)}+ \underset{N\to\infty}{O_\epsilon}(1) \ ,
 \end{align}
since $\displaystyle \left| \frac{\sin \big[ \Psi_N(y)-\Psi_N(x) \big]}{\pi(x-y)} \right| \le  2N + \epsilon^{-1/2}$ according to the estimate (\ref{Psi_2}). \\

It remains to show that the other terms in the  expansion (\ref{gue_1}) are uniformly bounded in $N$ for all $|x|\le1-\epsilon/\sqrt{2}$. By the triangle inequality,
\begin{equation*}
\big|\tilde\Lambda_N(x)\Lambda_N(y)-\tilde\Lambda_N(y)\Lambda_N(x) \big|=
\big|\tilde\Lambda_N(x)\big|\big|\Lambda_N(y)-\Lambda_N(x)\big| + \big|\tilde\Lambda_N(x)- \tilde\Lambda_N(y)\big|\big|\Lambda_N(x)\big| \ .
\end{equation*}
By (\ref{Lambda}), we see that $\big|\Lambda_N(y)-\Lambda_N(x)\big| \le C \epsilon^{-5/2} |x-y| $ 
and this estimates holds for $\tilde\Lambda_N$ as well, so that
\begin{equation}\label{gue_2}
\left|\frac{\tilde\Lambda_N(x)\Lambda_N(y)-\tilde\Lambda_N(y)\Lambda_N(x)}{x-y} \right| \le 2C^2N^{-1} \epsilon^{-4} \ .
\end{equation}
We can handle the third term similarly, by the triangle inequality and (\ref{Lambda}),
\begin{align*}
\big|\cos \big[\Psi_N(x)\big] \Lambda_N(y) - \cos \big[\Psi_N(y)\big] \Lambda_N(x) \big|
&\le \big| \Lambda_N(y) - \Lambda_N(x)\big| +   \big| \cos\Psi_N(x)- \cos\Psi_N(y)\big|  \big|\Lambda_N(y)  \big| \\
&\le C \epsilon^{-5/2} |x-y|  + C N^{-1} \epsilon^{-3/2}  \big|\Psi_N(x)- \Psi_N(y)\big|  
\end{align*}
Then, by (\ref{Psi_2}), we obtain the upper-bound,
\begin{equation}\label{gue_3}
\left|\frac{\cos \big[\Psi_N(x)\big] \Lambda_N(y) - \cos \big[\Psi_N(y)\big] \Lambda_N(x)}{x-y} \right| \le  4C \epsilon^{-5/2} \ .
\end{equation}
The fourth term  is also uniformly bounded since (\ref{gue_3}) holds for $\tilde\Psi_N$ and $\tilde\Lambda_N$ as well. Hence, if we put together $(\ref{gue_1} - \ref{gue_3})$, we have proved that
\begin{equation*} 
\sqrt{2N}K_N^{\omega_G}(\sqrt{2N}x, \sqrt{2N} y)  
= \frac{ \sin \big[ \Psi_N(y)-\Psi_N(x) \big]}{\pi(x-y)}+ \underset{N\to\infty}{O_\epsilon}\left(1 + \frac{1}{|x-y|N^2}\right) \ .
\end{equation*}

By (\ref{H}) and (\ref{H_1}), for all $|x| <1$, 
\begin{align*}
 \Psi_N(y)-\Psi_N(x)  = 2N \int_y^x \sqrt{1-t^2} dt + \frac{\arccos y - \arccos x}{2} \ ,
\end{align*}
and
\begin{equation*} 
\sqrt{2N}K_N^{\omega_G}(\sqrt{2N}x, \sqrt{2N} y)  
= \frac{ \sin \big[\displaystyle 2N \int_y^x \sqrt{1-t^2} dt \big]}{\pi(x-y)}+ \underset{N\to\infty}{O_\epsilon}\left(1 + \frac{1}{|x-y|N^2}\right) \ .
\end{equation*}
If we go back to the GUE correlation kernel, by formula (\ref{K_GUE}), we conclude that for all $|x|\le\sqrt{2}-\epsilon$, 
\begin{equation*} \label{gue_5}
K_N^{\G}(x,y) 
= \frac{ \sin \big[\displaystyle 2N \int_{y/\sqrt{2}}^{x/\sqrt{2}} \sqrt{1-t^2} dt \big]}{\pi(x-y)}+ \underset{N\to\infty}{O_\epsilon}\left(1 + \frac{1}{|x-y|N^2}\right) \ .
\end{equation*}
To obtain formula (\ref{gue_0}), it remains to make the change of variable $u=\sqrt{2} t$ in the last integral, then the semicircular law, (\ref{sc}), appears naturally. \qed\\

\subsection{The local asymptotics and uniformity} \label{sect:GUE_3}
 
The asymptotics of lemma~\ref{thm:sine_meso'}  is not uniform and to complete the proof of formula (\ref{sine_GUE}) below, we need to remove the condition $|x-y| \gg 1/N^2$. To do so, we will use a method introduced by Levin and Lubinsky to prove local universality, see~\cite{LL_09,Lubinsky_09b}. It consists in first computing the asymptotics of the Christoffel-Darboux kernel along the diagonal, then extending the result off-diagonal using some a priori estimates on the derivative of the OPs. For the GUE kernel, we can use that the Hermite function solves a second order ODE to obtain this estimate, see formula (\ref{Schrodinger}) and lemma~\ref{thm:sine_micro'} below.

\begin{proposition}\label{thm:sc} 
For any $|x|\le\sqrt{2}-\epsilon$,
\begin{equation}\label{sc_0}
 K_N^{\G}(x,x)= \frac{N\sqrt{2-x^2}}{\pi} + \underset{N\to\infty}{O}(1) \ .
 \end{equation}
\end{proposition}
 
 \proof  The Hermite polynomials are an Appell sequence and, by formula (\ref{Hermite_0'}), this implies that for all $k \ge 0$, 
 $$ \phi_k'(x) = \frac{k\gamma_k}{\gamma_{k-1}}  \phi_{k-1}(x)- x \phi_k(x) 
 = \sqrt{2k}  \phi_{k-1}(x)- x \phi_k(x)    \ .$$
 If we use this equation and formula (\ref{micro_2}) below, we obtain
 \begin{equation} \label{density_GUE}  
K_N^{\omega_G}(x,x) = N\left\{  \phi_{N-1}(x)^2- \sqrt{1-N^{-1}} \phi_{N}(x)\phi_{N-2}(x)  \right\} \ . 
\end{equation}
Let $\tilde \Psi_N(x) = \Psi_N(x) - \arccos x$. The same argument as the proof of proposition~\ref{thm:Hermite} shows that for any $|x|\le1-\epsilon/\sqrt{2}$,
\begin{equation*}\label{Hermite_1''}
\phi_{N-2}(\sqrt{2N}x)
=\frac{\eta_{N-1}}{\big(N(1-x^2) \big)^{1/4}}\bigg\{\cos \big[\tilde\Psi_N(x)-\arccos(x)\big]   
+\underset{N\to\infty}{O_\epsilon}\big(N^{-1}\big) \bigg\} \ .
\end{equation*}
By formulae (\ref{K_GUE})  and (\ref{density_GUE}), this implies that for any $|x|\le1-\epsilon/\sqrt{2}$,
\begin{align*}
&K_N^{\G}\left(\sqrt{2}x , \sqrt{2}x \right)=
\frac{N \eta_N^2}{\sqrt{1-x^2} }\bigg\{ \big(\cos \tilde\Psi_N(x)\big)^2\\
&\hspace{1.3cm}-  \sqrt{1-N^{-1}} \frac{\eta_{N+1}\eta_{N-1}}{\eta_N^2}\cos \big[\tilde\Psi_N(x)-\arccos(x)\big]\cos \big[\tilde\Psi_N(x)+\arccos(x)\big]
+O_\epsilon\big(N^{-1}\big) \bigg\}  
\end{align*}
By (\ref{nu}), 
$$  \eta_N^2 = \sqrt{2}/\pi + O(N^{-1}) \ ,
 \hspace{1cm}
  \sqrt{1-N^{-1}} \frac{\eta_{N+1}\eta_{N-1}}{\eta_N^2} = 1+ O(N^{-1}) \ ,$$ 
and using the trigonometric identity
$$ \cos \big[\tilde\Psi_N(x)-\arccos(x)\big]\cos \big[\tilde\Psi_N(x)+\arccos(x)\big]
= x^2-1+ \big(\cos \tilde\Psi_N(x)\big)^2 \ , $$
this yields for all $|x|\le1-\epsilon/\sqrt{2}$,
$$K_N^{\G}\left(\sqrt{2}x , \sqrt{2}x \right)
=\frac{N  \sqrt{2}}{\pi\sqrt{1-x^2} }\big\{ 1- x^2 + O_\epsilon\big(N^{-1}\big) \big\}  \ . $$
Formula (\ref{sc_0}) follows from a trivial change of variables. \qed\\

 \begin{lemma} \label{thm:sine_micro'}
 For any $\epsilon>0$, there exists two constants $A,C>0$ such that $|x| , |y|\le \sqrt{2} -\epsilon$,
\begin{equation} \label{sine_micro}
 \left| K_N^{\G}(x, y)- \frac{\sin\big[ N \pi \big( F_{\s}(x) - F_{\s}(y)\big)\big]}{\pi(x-y)} \right|
 \le A+ C |x-y| N^2 \ . 
 \end{equation}
 \end{lemma}
 
 \proof Let $I_N =[-\sqrt{N}(\sqrt{2}-\epsilon) , \sqrt{N}(\sqrt{2}-\epsilon) ]$.
 By  formula (\ref{K_G}),
\begin{equation} \label{micro_1} 
K_N^{\omega_G}(x, x+ \zeta) 
 = \sqrt{\frac{N}{2}}\frac{ \phi_{N-1}(x)\big(\phi_{N}(x+ \zeta)- \phi_{N}(x) \big) - \phi_{N}(x)\big(\phi_{N-1}(x+ \zeta)- \phi_{N-1}(x)\big)}{\zeta} \ ,
 \end{equation}
and taking the limit as $\zeta\to 0$, we get
\begin{equation} \label{micro_2}  
K_N^{\omega_G}(x,x) = \sqrt{N/2} \big\{  \phi_{N-1}(x)  \phi_{N}'(x) -  \phi_{N}(x)\phi_{N-1}'(x)  \big\} \ . 
\end{equation}
Then,  if we perform a $2^{nd}$ order Taylor expansion in formula (\ref{micro_1}), we get for all $x, x+\zeta\in I_N$, 
  \begin{equation}  \label{micro_3}     
 \big| K_N^{\omega_G}(x, x+ \zeta) - K_N^{\omega_G}(x,x) \big|  
 \le \sqrt{2N} |\zeta|
\max_{I_N} \big\{|\phi_N ''| , |\phi_{N-1} ''|  \big\}  \max_{I_N} \big\{|\phi_N | , |\phi_{N-1} |  \big\} \ .   
  \end{equation}
 The Hermite functions are known to solve the ODE:
\begin{equation} \label{Schrodinger}
\phi_k''(x) = (x^2-  2k+1) \phi_k(x) \  . 
\end{equation}
This implies that 
\begin{equation}  \label{micro_4}  
 \max_{I_N} \big\{|\phi_N ''| , |\phi_{N-1} ''|  \big\} \le (2N+1)  \max_{I_N} \big\{|\phi_N | , |\phi_{N-1} |  \big\} 
 \end{equation}
Moreover, by formula (\ref{Hermite_1}), there exists a constant $C$ which only depends on $\epsilon$ such that  
$$ \max_{I_N} \big\{|\phi_N | , |\phi_{N-1} |  \big\}   \le C/ 3N^{1/4} .$$ 
Hence, by (\ref{micro_3}) and  (\ref{micro_4}), we have  for all $x, x+\zeta \in I_N$,
  \begin{equation*}   
 \big| K_N^{\omega_G}(x, x+ \zeta) - K_N^{\omega_G}(x, x) \big| 
 \le C |\zeta| N \ .
 \end{equation*}
 By (\ref{K_GUE}), this implies that  for all $|x| , |x+\zeta |\le \sqrt{2} -\epsilon$,
  \begin{equation}  \label{micro_0}
 \big| K_N^{\G}(x, x+ \zeta) - K_N^{\G}(x, x) \big|   
 \le C |\zeta| N^2 \ ,
 \end{equation} 
and by formula (\ref{sc_0}) there exists another constant $A>0$ such that
  \begin{equation}   \label{micro_5}   
 \big|  K_N^{\G}(x, x+\zeta) - N \varrho_{\s}(x)  \big|  \le A +  C |\zeta| N^2 \ .
 \end{equation} 
On the other hand, by (\ref{sc}), the semicircular density is smooth and for all $|x| , |x+\zeta |\le \sqrt{2} -\epsilon$,
\begin{align*} \big|\sin\big[ N \pi \big( F_{\s}(x+\zeta) - F_{\s}(x)\big)\big] - \sin \big[ \pi N  \varrho_{\s}(x) \zeta \big]  \big| 
&\le \pi N \big| F_{\s}(x+\zeta) - F_{\s}(x) - \varrho_{\s}(x) \zeta |   \\
& \le \pi  |\zeta |^2 N \ .
\end{align*}
Moreover, if we use the trivial bound $| v-\sin v| \le v^2$,  
$$ \left|  \frac{\sin \big[ \pi N  \varrho_{\s}(x) \zeta \big] }{\pi\zeta} - N \varrho_{\s}(x) \right|  \le  |\zeta| N^2 \ ,$$
and by the triangle inequality, we get   for all $|x| , |x+\zeta |\le \sqrt{2} -\epsilon$,
\begin{equation} \label{micro_6}
 \left| \frac{\sin\big[ N \pi \big( F_{\s}(x+\zeta) - F_{\s}(x)\big)\big]}{\pi\zeta} - N \varrho_{\s}(x) \right| 
 \le 2|\zeta| N^2  \ .
\end{equation}
The lemma follows by combining the estimates (\ref{micro_5}) and (\ref{micro_6}). \qed\\
 
  \begin{remark} \label{rk:Levin_Lubinsky}  Using the same argument, it is possible to get the estimate $(\ref{sine_micro})$ for a general ensemble $\mathbb{P}_N^\omega$  provided that its correlation kernel correlation satisfies $(\ref{micro_0})$ and  
 $K^\omega_N(x,x)/N = \varrho_\omega(x) + O(N^{-1})$.  For instance,  if the weight $\omega$  do not depend on the dimension $N$ and is compactly supported, the estimate  $(\ref{micro_0})$ follows from the Markov-Berstein inequality, see \cite{Lubinsky_09b}. For the modified Jacobi ensembles, in the regime $\alpha>1/2$, this can be used to give another proof of proposition~$\ref{thm:sine_J}$ without using the local asymptotics $(\ref{sine_J_1})$, though it only gives an error term of order $N^{1/2-\alpha}$.
 \end{remark}

 By combining lemmas~\ref{thm:sine_meso'} and~\ref{thm:sine_micro'}, we obtain the full asymptotics for the GUE kernel. Notice that the error term of order $N^{-2}$ is crucial to complete the proof for all mesoscopic scale $\alpha\in(0,1]$. This can be achieved because  the asymptotics of proposition~\ref{thm:Hermite} includes an extra term compared to the classical expansion used in section~\ref{sect:UIE}.

  \begin{theorem} \label{thm:sine_GUE}
For any $L,\epsilon>0$ and for any $\alpha\in(0,1]$,
 \begin{equation} \label{sine_GUE}
\frac{1}{N^\alpha} K_N^{\G}\left(x_0 + \frac{\xi}{N^\alpha} , x_0 + \frac{\zeta}{N^\alpha} \right)= \frac{\sin \pi N \big( ( F_{\s} (x_0 + \xi N^{-\alpha} ) - F_{\s}(x_0 + \zeta N^{-\alpha} )\big)}{\pi(\xi-\zeta)} + \underset{N\to\infty}{O}(N^{-\alpha}) \ , 
 \end{equation}
uniformly for all $|x_0| \le \sqrt{2}-\epsilon$ and all $\xi,\zeta \in[-L,L]$.\end{theorem}
 
 \proof  We let $x=x_0 + \xi  N^{-\alpha}$ and $y=x_0 + \zeta  N^{-\alpha}$. In particular $|x-y| = N^{-\alpha}|\xi-\zeta|$ and, if $N$ is sufficiently large compared to $L$, then $|x|, |y| \le \sqrt{2}-\epsilon/2$.\\
By lemma~\ref{thm:sine_meso'}, if  $|\xi-\zeta| \ge N^{-2+\alpha}$,  we get
  \begin{equation*} 
\frac{1}{N^\alpha} K_N^{\G}\left(x,y\right)= \frac{\sin \pi N \big( ( F_{\s} (x_0 + \xi N^{-\alpha} ) - F_{\s}(x_0 + \zeta N^{-\alpha} )\big)}{\pi(\xi-\zeta)} + \underset{N\to\infty}{O}(N^{-\alpha}) \ .
 \end{equation*}
 On the other hand, if $|\xi-\zeta| \le B N^{-2+\alpha}$, then by lemma~\ref{thm:sine_micro'}, 
 \begin{equation*} 
 \left| \frac{1}{N^\alpha} K_N^{\G}(x, y)-\frac{\sin \pi N \big( ( F_{\s} (x_0 + \xi N^{-\alpha} ) - F_{\s}(x_0 + \zeta N^{-\alpha} )\big)}{\pi(\xi-\zeta)} \right|
 \le  \frac{A+BC}{N^\alpha}  \ . 
 \end{equation*}
 Since the parameter $B$ is arbitrary, we conclude that the error in formula (\ref{sine_GUE}) is uniform for all $\xi,\zeta \in [-L,L]$.\qed\\

\section*{Acknowledgment}
The author thanks Kurt Johansson for many valuable discussions and for reading carefully the draft of this  article.

\appendix

\section{Variance estimate in the global regime} \label{A:variance}

In this section, we consider the unitary invariant ensemble $\mathbb{P}^V_N$  introduced in theorem~\ref{thm:UIE} and we assume that  there exists $\mathrm{B}>2$ and $\eta>0$ so that 
\begin{equation} \label{potential'} 
 V(x) \ge 2 (1+\eta)  \log|x| \ , \hspace{1cm} \forall |x|> \mathrm{B} \ .  
\end{equation}
We also suppose that the potential $V$ satisfies the {\bf one-cut} condition and $J_V=(-1,1)$. So, we can apply the results of section~\ref{sect:variance}. 
We say that a real-valued function $f$ belong to the space  $\mathscr{H}^{1/2}$ and we denote $f\in \mathscr{H}^{1/2}$ if $f\in L^\infty(\R)$ and
\begin{equation} \label{H}
  \iint\limits_{[-1,1]^2} \left| \frac{f(x)- f(y)}{x-y} \right|^2 dxdy <\infty \ . 
  \end{equation}

\begin{lemma}\label{lemma_2} Let $f \in \mathscr{H}^{1/2}$ with compact support. If $f$ satisfies the condition $(\operatorname{\ref{L}})$ as well, then $f\in H^{1/2}(\R)$.
\end{lemma}

\proof 
Suppose that $\supp(f) \subseteq[-\mathrm{A},\mathrm{A}]$ and let $\mathcal{K}= \big\{|x| \le \mathrm{A}, 1\le |y|\le \mathrm{A}+1 \big\}$ and 
$\mathcal{B}= [-\mathrm{A},\mathrm{A}] \times [\mathrm{A}+1,\infty)$. By symmetry, we have
$$ \| f\|_{H^{1/2}}^2 \le \iint\limits_{[-1,1]^2}   \left| \frac{f(x)- f(y)}{x-y} \right|^2 dxdy  + 
2\iint\limits_{\mathcal{K}}   \left| \frac{f(x)- f(y)}{x-y} \right|^2 dxdy 
+4\iint\limits_{\mathcal{B}}   \left| \frac{f(x)- f(y)}{x-y} \right|^2 dxdy \ .$$
By \eqref{H}, the first term is finite. Since $f$ satisfies the condition (\ref{L}), the second term is bounded by $4 \mathrm{A}^2\mathrm{L}$. By definition of the set $\mathcal{B}$,  the third term satisfies\begin{align*} \iint\limits_{\mathcal{B}}   \left| \frac{f(x)- f(y)}{x-y} \right|^2 dxdy 
&\le \iint\limits_{\mathcal{B}}   \left| \frac{f(x)}{y-\mathrm{A}} \right|^2 dxdy \\
& \le 2\mathrm{A} \|f\|_\infty^2 ,
\end{align*}
and we conclude that 
$ \| f\|_{H^{1/2}}^2 <\infty$.\qed\\

The aim of this appendix is to derive an estimate for the variance of global linear statistics valid for continuously differentiable test functions.

\begin{proposition} \label{thm:variance_global} 
 Let $V:\R\to\R$ be a real-analytic function which satisfies $(\ref{potential'})$ and such that $J_V=(-1,1)$.  We denote $\Xi_N h= \sum h(\lambda_k)$ where the sum is over the eigenvalues of a random matrix distributed according to  $\mathbb{P}^V_N$. Let $h\in C^1(\R)$ and suppose that there exists $Q, n>0$ so that $|h'(x)|\le Q |x|^{n}$ for all $|x| \ge 1$, then
$$ \limsup{N\to\infty} \Var\big[\Xi_N h\big] \le  16 \tilde\Sigma(h)^2 \ ,$$
where
\begin{equation} \label{variance_2'} 
\tilde{\Sigma}(f)^2= \frac{1}{\pi^2}  \iint\limits_{[-1,1]^2} \left| \frac{f(x)- f(y)}{x-y} \right|^2 \frac{dxdy}{\sqrt{1-x^2}\sqrt{1-y^2}}  \ .
 \end{equation} 
\end{proposition}

The proof is based on the result of proposition~\ref{thm:variance} and the exponential decay of the Christoffel-Darboux kernel outside of the bulk; see lemma~\ref{lemma_1} below. We suppose that $h\in C^1(\R)$ in order to simplify the proof, however this condition is not necessary.  In fact, by a simple modification of our method, it suffices to suppose that $h \in \mathscr{H}^{1/2}$ and there exists $Q>0$ and $n>0$ so that for all $|x|>1-\delta$,
 \begin{equation*} 
  \sup\left\{  \left| \frac{h(x)- h(y)}{x-y} \right|  : |y| \le |x| \right\}
\le Q |x|^n \ .
\end{equation*}

\begin{lemma}\label{lemma_1} Under the assumption of proposition~\ref{thm:variance_global}, if we also suppose  that $h(x)=0$ for all $|x|\le \mathrm{B}$, cf.~formula $(\operatorname{\ref{potential'}})$. Then, there exists $C >0$ so that  
$$ \Var\big[\Xi_N h\big]  \le C \mathrm{B}^{- \eta N} \ .  $$
\end{lemma}


\proof By \cite[formula 1.58]{Deift_al_99_b}, for any $\epsilon>0$, we have\begin{equation}  \label{decay_1}
 \big|\Phi_N(x)\big| \le \left(\frac{1}{2\sqrt{\pi}}\left|\frac{x+1}{x-1}\right|^{1/4} + \underset{N\to\infty}{O_\epsilon}\big(N^{-1}\big)\right)  e^{-N \mathfrak{H}_V(x)} \ , \hspace{.6cm}\forall  |x|>1+\epsilon\ ,
\end{equation}
where for all $x\in\R$,
\begin{equation*}
\mathfrak{H}_V(x)=  \frac{V(x)+\ell}{2} - \int \log|x-s|\varrho_V(s) ds 
\hspace{.6cm}\text{and}\hspace{.6cm} \ell \in \R. 
\end{equation*}

This function appears in the determination of the equilibrium density $\varrho_V$. In fact, $\varrho_V(x)dx$  is the unique minimizer of a weighted energy functional and it is  uniquely determined by the following Euler-Lagrange variational conditions:
$$ \begin{cases}
 \mathfrak{H}_V(x) = 0  &\forall x\in \overline{J_V} \\
 \mathfrak{H}_V(x) \ge 0 &\forall x\in \R\backslash\overline{J_V}   \ .
\end{cases}$$
Moreover, since $\supp(\varrho_V)=[-1,1]$, we have for all $|x|>1$ 
$$ \int \log|x-s|\varrho_V(s) \le \log(2|x|) \  . $$  
Hence, if the potential $V(x)$ satisfies the condition (\ref{potential'}), then  $\mathfrak{H}_V(x) \ge \eta \log|x| +\frac{\ell}{2}-\log 2$ for all $|x|>\mathrm{B}$. In fact, choosing a larger constant $\mathrm{B}$  if necessary, we can suppose that $\mathfrak{H}_V(x) \ge \frac{\eta \log|x|}{2}$.  
 By formula (\ref{decay_1}), this implies that there exists $\mathrm{C}>0$ so that for all $|x|>\mathrm{B}$,
\begin{equation}  \label{decay_3}
 \big|\Phi_N(x)\big| \le  \sqrt{\mathrm{C}/2}\  e^{-N  \eta \log|x| / 2} \ .
 \end{equation}
Using~\cite[formula 1.59]{Deift_al_99_b} instead,  we can show that the estimate (\ref{decay_3}) holds  for the function $\Phi_{N-1}$ as well. By formula (\ref{variance_3}), this implies that  for all $|x| \ge \mathrm{B}$,
\begin{equation*} 
\big|K^\omega_N(x,y)\big|^2 \le \mathrm{C}  \frac{\gamma_{N-1}}{\gamma_{N}}
\frac{\big| \Phi_{N-1}(y)\big|^2+\big| \Phi_{N}(y)\big|^2}{|x-y|^2} e^{-N\eta \log |x| } \ .
\end{equation*}
Hence, since $\| \Phi_N \|_{L^2}= \| \Phi_{N-1} \|_{L^2}=1$, we obtain  for all $|x| \ge \mathrm{B}$,
\begin{equation}  \label{decay_4}
  \int_\R \big|(x-y)K^\omega_N(x,y)\big|^2 dy   
 \le   \mathrm{C}  \frac{\gamma_{N-1}}{\gamma_{N}}e^{-N\eta\log |x| } \ .   
 \end{equation}

On the other hand, by assumptions, we have for all $|y|\le|x|$, 
\begin{align*}  \left| \frac{h(x)- h(y)}{x-y} \right| & \le \1_{|x| >  \mathrm{B}} \sup\big\{  h'(t) : |t| \le |x|  \big\}  \\
&\le Q |x|^n  \1_{|x| > \mathrm{B}} \ .
\end{align*}
According to formula (\ref{variance_1}) and (\ref{decay_4}), we obtain
\begin{align} \notag
 \Var\left[ \Xi_N h \right] &= \frac{1}{2} \iint \left|h(x)-h(y) \right|^2  \big|K^\omega_N(x,y)\big|^2 dxdy \\
&\notag
 \le \frac{Q^2}{2} \int_{\R\backslash[-\mathrm{B},\mathrm{B}]} |x|^{2n} \left( \int_\R \big|(x-y)K^\omega_N(x,y)\big|^2 dy  \right) dx \\
 &\label{decay_6}
 \le \mathrm{C}Q^2 \frac{\gamma_{N-1}}{\gamma_{N}} \int_{\mathrm{B}}^\infty x^{2n} e^{-N\eta\log(x)} dx \ .
\end{align}
Because of the asymptotics (\ref{semiclassical_3}), $\displaystyle C:= \mathrm{C} Q^2\sup_{N\in\N}\{ \frac{\gamma_{N-1}}{\gamma_N}\}\mathrm{B}^{2n}<\infty$ and the proof is complete. \qed\\

\noindent{\it Proof of proposition~\ref{thm:variance_global}.}  Let $\mathrm{A}>\mathrm{B}$  and  $\chi \in C^1\big(\R_+ \to [0,1]\big)$ such that  $ -\chi' \in[0,1]$ and
\begin{equation*} \chi(x) =
\begin{cases} 1 &\text{if}\ x \le  \mathrm{B} \\  0 &\text{if}\ x \ge \mathrm{A}
\end{cases} \ .
\end{equation*} 
We also let $\tilde{\mathrm{A}}=\mathrm{A}+1$. We decomposition $h= f+g$ where $f =\chi h \in C^1_0(\R)$ and $g=(1-\chi) h \in C^1(\R)$. According to formula (\ref{variance_1}), we have
\begin{equation}\label{decay_7}
  \Var\big[ \Xi_N h \big]  \le 2 \left(   \Var\big[ \Xi_N f \big] +  \Var\big[ \Xi_N g \big] \right) \ . 
 \end{equation}
First, since $g(x)=0$ for all $|x|\le \mathrm{B}$ and $|g'(x)|\le |h(x)| + |h'(x)|$, by assumptions there exists a constant $\tilde{Q}$ so that $|g'(x)|  \le \tilde{Q} |x|^{n+1}$ for all $|x|\ge 1$. Then, by lemma~\ref{lemma_1},
\begin{equation}\label{decay_8}
 \limsup{N\to\infty }\Var\big[\Xi_N g\big] =0 \ .  
 \end{equation}
Next, we will show that the function $f$ satisfies the condition  (\ref{L}). By definition, we have
\begin{equation*}   \left| \frac{f(x)- f(y)}{x-y} \right|  \le  \left| \frac{h(x)- h(y)}{x-y} \right| \chi(x) +  |h(y)| \left| \frac{\chi(x)- \chi(y)}{x-y} \right| \ .
\end{equation*}
Hence, if $|x|< \tilde{\mathrm{A}}$, using the properties of the cutoff  function $\chi$, for all $ |y| \le |x|$,
\begin{equation*}   \left| \frac{f(x)- f(y)}{x-y} \right| 
\le  \sup\big\{ |h'(t)|+|h(t)| : |t| \le \mathrm{A}' \big\} \ .
\end{equation*}
On the other hand, if $|x|\ge \tilde{\mathrm{A}}$, for all $ |y| \le |x|$,
$$   \left| \frac{f(x)- f(y)}{x-y} \right|  \le  \frac{  |h(y)|\chi(y)}{|x-y|}  \le 
\begin{cases} 0 &\text{if}\ |y| \ge \mathrm{A} \\
 \sup\big\{ |h(t)| : |t| \le \mathrm{A}\big\} &\text{else}
\end{cases} $$
Hence, there exists $\mathrm{L}>0$ so that
\begin{equation*} 
 \sup \left\{ \left| \frac{f(x)- f(y)}{x-y} \right|  :  |y| \le |x|  \right\} = \mathrm{L} \ ,
 \end{equation*}
and, by symmetry, the function $f$ satisfies the condition (\ref{L}).
Moreover, by lemma~\ref{lemma_2}, $f\in H^{1/2}(\R)$ and  by proposition~\ref{thm:variance}, this implies that 
\begin{equation} \label{decay_9}
 \limsup{N\to\infty} \Var\left[ \Xi_N f \right]  \le 8\tilde\Sigma(f)^2 \ .
\end{equation}
Combining the estimates (\ref{decay_7} - \ref{decay_9}), we conclude that 
\begin{equation*}
 \limsup{N\to\infty} \Var\left[ \Xi_N h \right]  \le 16\tilde\Sigma(f)^2 \ .
\end{equation*}
It completes the proof since $\tilde\Sigma(f)=\tilde\Sigma(h)$  because 
 $h(x)=f(x)$ for all $|x| \le 1$.    \qed\\

Proposition~\ref{thm:variance_global} is used in \cite{L_15b} to give a new proof of  theorem~\ref{thm:clt_2} below.
 In fact, the results of \cite{L_15b} are valid for more general orthogonal polynomial ensembles. 
 Theorem~\ref{thm:clt_2} is an extension of the CLT (\ref{clt_2}) and  its proof is inspired from that of theorem~\ref{thm:clt_1}.

\begin{theorem}\label{thm:clt_2}  Let $V:\R\to\R$ be a real-analytic function which satisfies the condition $(\ref{potential'})$ and such that $J_V=(-1,1)$.  
 If $(\lambda_1, \dots, \lambda_N)$ denote  the eigenvalues of a random matrix distributed according to  $\mathbb{P}^V_N$, then for any $f\in C^1(\R)$ such that there exists $Q,n>0$ so that $|f'(x)|\le Q |x|^{n}$ for all $|x| \ge 1$, we have
\begin{equation} \label{clt_2'}
\sum_{k=1}^N f(\lambda_k)  - \E{\sum_{k=1}^N f(\lambda_k)}  \underset{N\to\infty}{\Longrightarrow}   \No\big( 0, \Sigma(f)^2 \big)    \ .
\end{equation}
\end{theorem}

\proof For any $f\in C^1(\R)$, we denote  $\displaystyle\overline{\Xi_N}f = \Xi_Nf - \E{ \Xi_Nf }$ and 
$\vartheta_N(\xi; f)= \E{e^{i \xi \overline{\Xi_N} f }}$.
Proposition~\ref{thm:variance_global} implies that the sequence of random variables  $\overline{\Xi_N} h$ is tight and by Prokhorov's theorem, there exists an increasing map $\pi:\N\to\N$ and a random variable  we denote $\mathfrak{S}(h)$  so that  $\overline{\Xi_{\pi(N)}}h \Rightarrow \mathfrak{S}(h)$ as $N\to\infty$. For any $\epsilon>0$, by Weierstrass's approximation theorem, there exists a polynomial $P_\epsilon$ so that
$$ \sup\big\{ | h'(x)-P'_\epsilon(x)| : |x| \le 1 \big\} \le \sqrt{\epsilon} \ , $$  
and by formula (\ref{variance_2'}),
\begin{equation*}
 \tilde\Sigma(h-P_\epsilon) \le \frac{\epsilon}{\pi^2}  \iint\limits_{[-1,1]^2} \frac{dxdy}{\sqrt{1-x^2}\sqrt{1-y^2}} = \epsilon \ .  
 \end{equation*}
Let $\chi \in C^\infty\big(\R_+ \to [0,1]\big)$ such that  
$\displaystyle \chi(x) =
\begin{cases} 1 &\text{if}\ x \le 1 \\  0 &\text{if}\ x \ge 2
\end{cases}$. For any $\epsilon>0$, the function $H_\epsilon=\chi P_\epsilon \in  C^\infty \cap L^\infty(\R)$ and by (\ref{clt_2}),
$$ \lim_{N\to\infty} \vartheta_N( \xi ; H_\epsilon) = \vartheta( \xi ; H_\epsilon) \ , $$
where  $\vartheta( \xi ; f)= e^{- \xi^2 \Sigma(f)^2/2}$  for any $f\in C^1(\R)$ and $\xi\in\R$. Moreover,  since $H_\epsilon=P_\epsilon$ on $[-1,1]$, $\tilde\Sigma(h-H_\epsilon)= \tilde\Sigma(h-P_\epsilon)$ and
\begin{equation}\label{tightness_1}
 \tilde\Sigma(h-H_\epsilon) \le \epsilon \ .
 \end{equation}
 By definition of the random variable $\mathfrak{S}(h)$, 
\begin{equation*}
 \left| \E{e^{i \xi \mathfrak{S}(h)}}-  \vartheta(\xi; H_\epsilon) \right|  =  \lim_{N\to\infty} \big| \vartheta_{\pi(N)}( \xi ; h) - \vartheta_{\pi(N)}( \xi ; H_\epsilon) \big|  \ ,
 \end{equation*}
and using the estimate (\ref{tightness}), we obtain
\begin{equation} 
 \left| \E{e^{i \xi \mathfrak{S}(h)}}-  \vartheta(\xi; H_\epsilon) \right|  
 \le  4|\xi| \limsup{N\to\infty} \sqrt{\Var\big[\Xi_{\pi(N)}h - \Xi_{\pi(N)}H_\epsilon\big]} \ .
 \end{equation} 
   Since the processes $\Xi_N$ are linear and the function $h-H_\epsilon$ satisfies the hypothesis of proposition~\ref{thm:variance_global}, for any subsequence $\pi$,
\begin{equation}\label{tightness_2}
 \limsup{N\to\infty} \sqrt{\Var\big[\Xi_{\pi(N)}h - \Xi_{\pi(N)}H_\epsilon\big]} \le  16  \tilde\Sigma(h-H_\epsilon) \ .
 \end{equation}
By formulae (\ref{tightness_1} - \ref{tightness_2}), this implies that
\begin{equation} \label{tightness_3}
 \left| \E{e^{i \xi \mathfrak{S}(h)}}-  \vartheta(\xi; H_\epsilon) \right|  
  \le 64|\xi| \epsilon \ . 
 \end{equation}
It is easy to check that $\Sigma$ and $\tilde\Sigma$ are  semi-norms (modulo constant) on the space $C^1([-1,1])$ and that  $\Sigma \le \tilde\Sigma$. Hence, by the triangle inequality,
$$
 \big|\Sigma(H_\epsilon)-  \Sigma(h)\big| \le  \tilde\Sigma(h-H_\epsilon) 
$$
By  (\ref{tightness_1}), it implies that $ \Sigma(H_\epsilon)\to \Sigma(h)$ as $\epsilon \to 0$ and for any $\xi\in\R$,
\begin{equation} \label{tightness_4}
 \lim_{\epsilon\to0}   \vartheta(\xi; H_\epsilon)= \vartheta(\xi; h) \ .
 \end{equation}
Combining (\ref{tightness_3}) and  (\ref{tightness_4}), we conclude that $\mathfrak{S}(h) \sim  \No\big(0,\Sigma( h)\big)$ and the CLT follows since this holds for any subsequence $\pi$. \qed\\

\begin{remark}\label{rk:GUE} For the GUE kernel, using Cram\'{e}r's inequality, $\|\phi_k\|_\infty \le 1$ for all $k\ge 0$, so that
 $$\left| K_N^{\G}(x,y) \right| \le N \ , \hspace{.8cm} \forall x,y \in\R \ . $$
Moreover, by Theorem~5.2.3 in \cite{Pastur_Shcherbina},  for any $\epsilon>0$ there exists $\beta, C>0$ so that $K_N^{\G}(x,x)\le C N e^{- \beta N x^2}$ for all $|x|\ge 1+\epsilon$. Hence, by the Cauchy-Schwartz inequality, we obtain for all $|x| \ge 1+\epsilon$ and $y\in\R$,
\begin{align*} \left| K_N^{\G}(x,y) \right|^2 & \le K_N^{\G}(y,y)  K_N^{\G}(x,x) \\ 
& \le C N^2 e^{- \beta N x^2}  \ .
\end{align*}
This implies that proposition~\ref{thm:variance_global} and the CLT (\ref{clt_2'}) hold for any test function  $h(x)= \underset{x\to\infty}{o}\big(e^{ |x|^\alpha} \big)$ with $0<\alpha<2$ and  such that there exists $0<\delta<1$ and $\mathrm{L}>0$ so that 
 \begin{equation*} 
  \sup\left\{  \left| \frac{h(x)- h(y)}{x-y} \right|  : |y| \le |x| ,\ 1-\delta <|x| <1+\delta  \right\}
\le \mathrm{L}\ .
\end{equation*}
\end{remark}

\bibliographystyle{siam}


\end{document}